\documentclass[12pt]{amsart}

\usepackage{amssymb,amscd,psfig,latexsym,graphicx,eucal}
\usepackage{fullpage}
\newtheorem{theorem}{Theorem}[section]
\newtheorem{lemma}[theorem]{Lemma}

\newtheorem{corollary}[theorem]{Corollary}

\theoremstyle{definition}
\newtheorem{definition}[theorem]{Definition}

\newtheorem{remark}[theorem]{Remark}

\numberwithin{equation}{section}

\newcommand{\til}{\triangleleft}
\newcommand{\hrp}{\widehat{R}_M}
\newcommand{\rp}{R_M}
\newcommand{\ds}{\displaystyle}
\newcommand{\rd}{R_c}
\newcommand{\per}{\operatorname{per}}
\newcommand{\Hdim}{\dim_{\operatorname{H}}}
\newcommand{\db}{{\bf d}}

\newcommand{\ve}{\varepsilon}
\newcommand{\es}{\emptyset}
\newcommand{\sm}{\smallsetminus}
\newcommand{\bd}{\partial}

\newcommand{\ov}{\overline}
\newcommand{\io}{\iota}
\newcommand{\con}{\operatorname{const.}}
\newcommand{\CC}{{\mathbb C}}

\newcommand{\RR}{{\mathbb R}}

\newcommand{\TT}{{\mathbb T}}
\newcommand{\ZZ}{{\mathbb Z}}

\newcommand{\DD}{{\mathbb D}}
\newcommand{\QQ}{{\mathbb Q}}
\newcommand{\vs}{\vspace{2mm}}

\newcommand{\iso}{\stackrel{\cong}{\longrightarrow}}

\renewcommand{\marginpar}[1]{}
\catcode`\@=12

\def\Empty{}
\newcommand\oplabel[1]{
  \def\OpArg{#1} \ifx \OpArg\Empty {} \else
  	\label{#1}
  \fi}
		
\long\def\realfig#1#2#3#4{
\begin{figure}[tp]
\centerline{\psfig{figure=#2,width=#4}}
\caption[#1]{#3}
\oplabel{#1}
\end{figure}}

\newcommand{\comm}[1]{}

\input{psbox}

\psfordvips

\newenvironment{pf}{\begin{proof}}{\end{proof}}
\newcommand{\cal}{\mathcal}
\newcommand{\thmref}[1]{Theorem~\ref{#1}}

\newcommand{\lemref}[1]{Lemma~\ref{#1}}

\newcommand{\corref}[1]{Corollary~\ref{#1}}
\newcommand{\figref}[1]{Fig.~\ref{#1}}

\begin{document}

\title[The real slice of the Mandelbrot set]
{External rays and the real slice of the Mandelbrot set}

\author[S. Zakeri]
{Saeed Zakeri}

\address{S. Zakeri, Department of Mathematics, University of
Pennsylvania, Philadelphia, PA 19104-6395, USA}

\curraddr{Institute for Mathematical Sciences, Stony Brook
University, Stony Brook, NY 11790-3651, USA}

\email{zakeri@math.sunysb.edu}



\date{August 20, 2002}

\begin{abstract}
This paper investigates the set of angles of the parameter rays which
land on the real slice $[-2,1/4]$ of the Mandelbrot set. We prove
that this set has zero length but Hausdorff dimension 1. We
obtain the corresponding results for the tuned images of the real
slice. Applications of these estimates in the study of critically
non-recurrent real quadratics as well as biaccessible points of
quadratic Julia sets are given.
\end{abstract}

\maketitle
\thispagestyle{empty} \def\IMSmarkvadjust{0 pt}
\def\IMSmarkhadjust{0 pt}
\def\IMSmarkhpadding{0 pt}
\def\IMSpubltext{Published in modified form:}
\def\SBIMSMark#1#2#3{
 \font\SBF=cmss10 at 10 true pt
 \font\SBI=cmssi10 at 10 true pt
 \setbox0=\hbox{\SBF \hbox to \IMSmarkhpadding{\relax}
                Stony Brook IMS Preprint \##1}
 \setbox2=\hbox to \wd0{\hfil \SBI #2}
 \setbox4=\hbox to \wd0{\hfil \SBI #3}
 \setbox6=\hbox to \wd0{\hss
             \vbox{\hsize=\wd0 \parskip=0pt \baselineskip=10 true pt
                   \copy0 \break%
                   \copy2 \break% 
                   \copy4 \break}}
 \dimen0=\ht6   \advance\dimen0 by \vsize \advance\dimen0 by 8 true pt
                \advance\dimen0 by -\pagetotal
	        \advance\dimen0 by \IMSmarkvadjust
 \dimen2=\hsize \advance\dimen2 by .25 true in
	        \advance\dimen2 by \IMSmarkhadjust

%
%   Check for publication info
%
%  \newread\jref
  \openin2=publishd.tex
  \ifeof2\setbox0=\hbox to 0pt{}
  \else 
     \setbox0=\hbox to 3.1 true in{
                \vbox to \ht6{\hsize=3 true in \parskip=0pt  \noindent  
                {\SBI \IMSpubltext}\hfil\break
                \input publishd.tex 
                \vfill}}
  \fi
  \closein2
  \ht0=0pt \dp0=0pt
 \ht6=0pt \dp6=0pt
 \setbox8=\vbox to \dimen0{\vfill \hbox to \dimen2{\copy0 \hss \copy6}}
 \ht8=0pt \dp8=0pt \wd8=0pt
 \copy8
 \message{*** Stony Brook IMS Preprint #1, #2. #3 ***}
}
 
\def\IMSmarkvadjust{-30pt}
\SBIMSMark{2002/02}{August 2002}{}

\tableofcontents

\section{Introduction}
\label{sec:intro}

The {\it Mandelbrot set} $M$ is the connectedness locus of the
family $Q_c:z \mapsto z^2+c$ of normalized complex quadratic polynomials:
$$M:=\{ c \in \CC : \operatorname{The\ Julia\ set\ of\ } Q_c \
\operatorname{is\ connected} \}.$$
It is a compact, full, and
connected subset of the plane, with a tremendously intricate
structure near the boundary (see \figref{M}). In recent years,
great deal of research has gone towards understanding the topology,
geometry, and combinatorics of $M$, as $M$ and its higher degree
cousins are the universal objects which appear in the bifurcation
locus of any holomorphic family of rational maps \cite{McMullen}.

The normalized Riemann mapping $\Phi : \widehat{\CC} \sm M \iso
\widehat{\CC} \sm \ov{\DD}$ which satisfies $\Phi (\infty)=\infty$ and
$\Phi'(\infty)=1$ plays a special role in the study of the
quadratic family and has a dynamical meaning: $\Phi (c)$ is the
conformal position of the escaping critical value $c$ in the basin
of attraction of infinity for $Q_c$. The normalized Lebesgue
measure $m$ on the circle $\TT=\RR / \ZZ \cong \bd \DD$ pulls back by
$\Phi$ to the {\it harmonic measure} $\mu_M$
supported in $\bd M$. More precisely, define the {\it parameter
ray} $\rp(t)$ of {\it external angle} $t \in \TT$ as the
$\Phi$-preimage of the radial line $\{ re^{2\pi i t}: r>1 \}$. We
say that $\rp(t)$ {\it lands} at $c \in \bd M$ if $\lim_{r \to 1}
\Phi^{-1}(re^{2\pi i t})=c$. It follows from a classical theorem of
\mbox{A.~Beurling} (see for example \cite{Pommerenke}) that every
parameter ray lands at a well-defined point of $\bd M$, except
possibly for a set of external angles of capacity zero
(conjecturally empty in this case). For a Borel measurable set $S \subset
M$, the harmonic measure is given by
%%%%%%%%%%%%%%%%%%%%%%%%
\realfig{M}{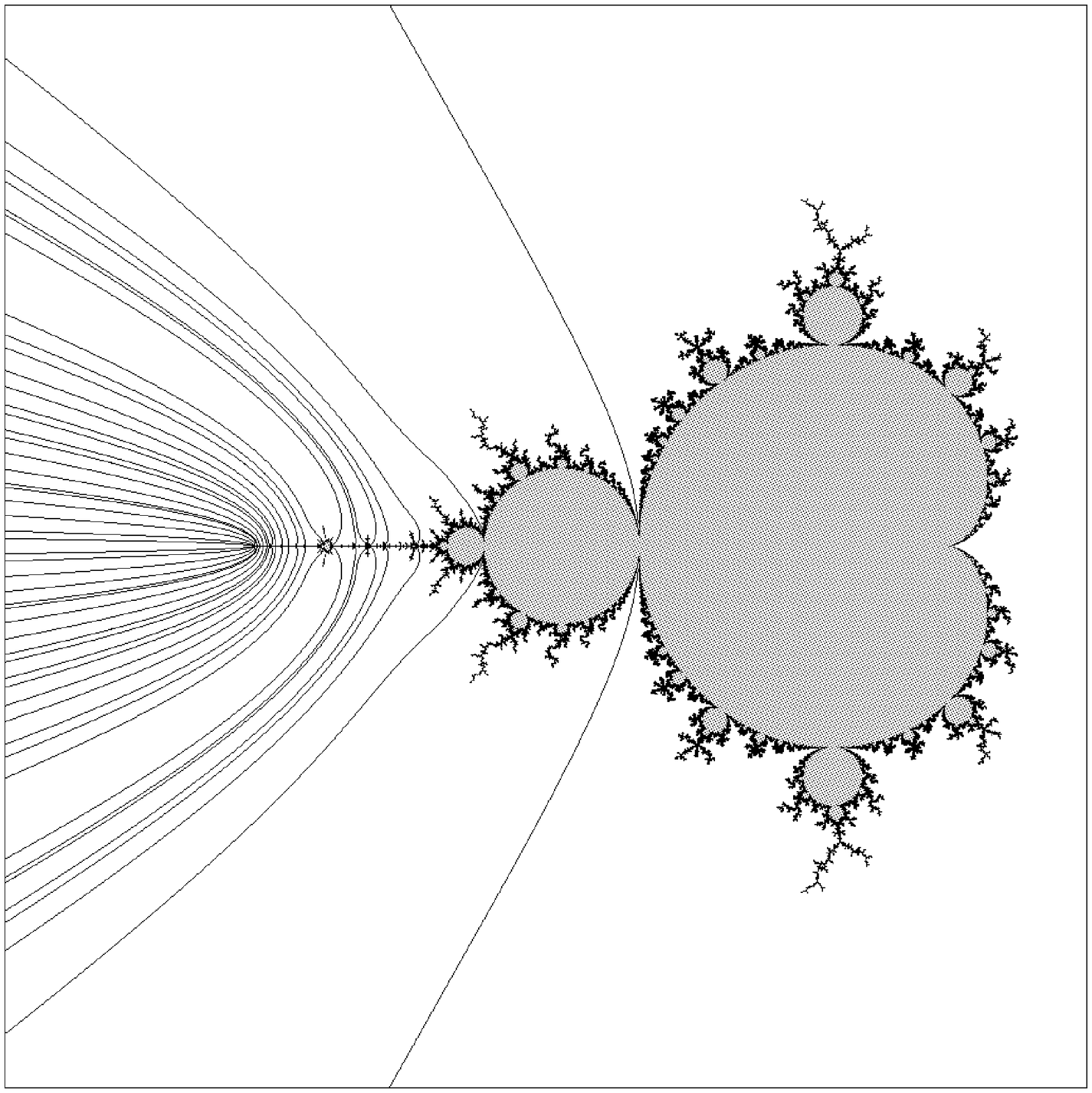}{\sl The Mandelbrot set M and some parameter rays
landing on its real slice.}{9cm}
%%%%%%%%%%%%%%%%%%%%%%%%
$$\mu_M(S) := m \{ t \in \TT : \rp(t)\ \operatorname{lands\ at\ a\
point\ of} S \} .$$
In other words, $\mu_M(S)$ is the probability
that a Poincar\'e geodesic in $\widehat{\CC} \sm M$ emanated from
infinity in a random direction hits $S$. The great complexity of
$M$ near each of its boundary points suggests that the harmonic
measure of embedded arcs in $M$ is zero so that they are almost
invisible from infinity. For example, all the hyperbolic components
of $M$ (the ``bulbs'' in \figref{M}) have piecewise analytic
boundary and a direct computation shows that the harmonic measure
of each of these boundary arcs is zero \cite{Douady1}. On the other hand,
$M$ contains many other essentially different embedded arcs, the
basic example of which is the {\it real slice} $M \cap
\RR=[-2,1/4]$. By tuning (see \S \ref{sec:bg} below), one obtains
countably many distinct embedded arcs in $M$ which are all images
of the real slice under embeddings $M \hookrightarrow M$.
The harmonic measure of each tuned image of the real slice is zero
since it can be shown using the tuning algorithm that the
corresponding angles are contained in a self-similar Cantor set
(see \cite{Manning} and compare \S \ref{sec:bg}). However, this
question for the real slice itself is non-trivial:

\begin{theorem}\label{A}
The harmonic measure of the real slice of the Mandelbrot
set is zero.
\end{theorem}

Naturally, one is led to consider a finer tool to measure the size of the
real slice as seen from infinity. In contrast to \thmref{A}, we show that

\begin{theorem}\label{B}
The set of external angles $t\in \TT$ for which the parameter ray
$\rp(t)$ lands at a point of $[-2, 1/4]$ has Hausdorff dimension $1$.
\end{theorem}

The main ingredient of the proofs of these theorems is an explicit
description for the set $\cal R$ of all angles $t \in [0,1/2]$ such
that the prime-end impression of the ray $\rp(t)$ in $\bd M$
intersects the real line (see \S \ref{sec:quad}). Conjecturally,
this set coincides with the set of angles of the parameter rays
which land on the real slice, but whether or not this is true is
irrelevant here, as the difference between the two sets has zero
capacity by Beurling's theorem. The description of $\cal R$ is
obtained by establishing the existence of a canonical homeomorphism
$c \mapsto \tau(c)$ between $\bd M \cap {\RR}$ and $\cal R$, with
the property that the prime-end impression of the parameter ray at
angle $\tau(c)$ intersects $\RR$ precisely at $c \in \bd M$;
alternatively, the associated dynamic ray at angle $\tau(c)$ lands
at the critical value $c$ in the Julia set of $Q_c$ (see
\thmref{2rays} and \lemref{R=c}).

Once the appropriate
description of $\cal R$ is in hand, \thmref{A} follows easily from
ergodicity of the doubling map $\db:t \mapsto 2t$ (mod 1) on the
circle. For the Hausdorff dimension question, we introduce a
one-parameter family of compact sets $\{ {\cal K}_{\sigma} \}_{\sigma >0}$
in \S \ref{sec:AB} whose dimension is estimated from
below by an application of Frostman's Lemma. The close
relation between the family $\{ {\cal K}_{\sigma} \}_{\sigma >0}$
and the set $\cal R$ allows
us to use these estimates and prove \thmref{B}.

The discussion in \S \ref{sec:AB} concludes with a generalization of
\thmref{B} to all tuned images of the real slice:

\begin{theorem} \label{C}
Let $H$ be a hyperbolic component of $M$ of period $p>1$ and
let $\eta_H \subset M$ be the corresponding tuned image of the real slice
$[-2,1/4]$. Then, the set of external angles $t \in \TT$  for which
$\rp(t)$ lands at a point of $\eta_H$ has Hausdorff dimension
$1/p$.
\end{theorem}

Combining \thmref{B} with standard dimension theorems in conformal
mapping theory leads to dimension estimates in the parameter space
of real quadratic polynomials. As an example, we show in \S
\ref{sec:gen} that

\begin{theorem} \label{CNR}
The set of parameters $c \in \bd M \cap \RR$ for which the quadratic
$Q_c$ is critically non-recurrent has Hausdorff dimension $1$.
\end{theorem}

Note that this set has Lebesgue measure
zero by a theorem of D.~Sands, although the full set $\bd M \cap \RR$
has positive Lebesgue measure according to Jakobson (compare \cite{Sands}
and \cite{Jakobson}).

Other applications of these estimates will be discussed in \S
\ref{sec:bi}. Let $B_c$ denote the set of angles of dynamic rays
which land on the {\it biaccessible} points in the Julia set
of the quadratic polynomial $Q_c$. In other words, $t \in
B_c$ if there exists an $s \neq t$ such that the dynamic rays at
angles $t$ and $s$ land at a common point of the Julia set.

\begin{theorem} \label{D}
For $-2 < c \leq -1.75$,
$$0< \ell(c) \leq \Hdim(B_c) < 1,$$
where $\ell(c)$ is an explicit constant which tends to $1$ as
$c$ tends to $-2$. In particular,
$$\lim_{c\,{\scriptstyle\searrow}\, -2} \Hdim (B_c) = \Hdim (B_{-2})= 1.$$
\end{theorem}

The function $c \mapsto \Hdim(B_c)$ is monotonically decreasing on
$[-2,1/4]$ and vanishes for $c_{\operatorname{Feig}} < c \leq 1/4$, where
$c_{\operatorname{Feig}} \approx -1.401155$ is the Feigenbaum
value (see \S \ref{sec:bi}). I do not know what exactly happens for
$-1.75 < c \le c_{\operatorname{Feig}}$.

The explicit form of $\ell(c)$ in \thmref{D}
will be given in \S \ref{sec:bi}. It
is interesting to contrast this theorem  with the fact that the
measure of $B_c$ is zero for all complex parameters $c \neq -2$ (see
\cite{Stas}, \cite{Zakeri}, \cite{Zdunik}). The statement that
$B_c$ has positive Hausdorff dimension has been shown by S.~Smirnov
for Collet-Eckmann real quadratics by a very different argument
\cite{Stas}.\\ \\
{\bf Acknowledgment.} I would like to thank \mbox{Peter Jones} who,
some time ago, suggested the study of the parameter rays landing on
the real slice. I am grateful to Stas Smirnov who brought
\cite{Sands} to my attention and kindly pointed out that a more
careful application of the argument in my first draft gives the
sharper result of \thmref{CNR}. My further thanks are due to John
Milnor and Carsten L.~Petersen for their useful comments on earlier
drafts of this paper. The pictures of the Julia sets and the
Mandelbrot set with its parameter rays are created by John Milnor's
program {\tt polyjul} and Scott Sutherland's adapted version {\tt
manray}.

\section{Background material}
\label{sec:bg}

We collect a few basic facts about quadratic Julia sets,
hyperbolic components of the Mandelbrot set,
rational parameter rays, and the tuning algorithm.
For details, see \cite{Douady2}, \cite{Douady-Hubbard},
\cite{Haissinsky}, and \cite{Milnor}.

\subsection*{Quadratic Julia sets}
\label{subsec:qjs}

Fix a parameter $c \in \CC$ and consider the quadratic polynomial
$Q_c: z \mapsto z^2+c$. The {\it filled Julia set} of $Q_c$,
denoted by $K_c$, is the set of all points in the plane with bounded
forward orbit under $Q_c$. The topological
boundary $J_c:=\bd K_c$ is called the {\it
Julia set} of $Q_c$. Both sets are non-empty, compact, and
totally-invariant.
Moreover, $K_c$ is always full, in the sense that $\CC \sm K_c$ is
connected. The domain $\CC \sm K_c$ is called the
{\it basin of attraction of infinity}; it consists of all
points with forward orbit tending to $\infty$. The components of the
interior of $K_c$, if any, are called the {\it bounded Fatou components}
of $Q_c$.

When $c \in M$, the filled Julia set $K_c$ is connected,
so there exists a unique conformal isomorphism
$\varphi_c: \widehat{\CC} \sm K_c \iso \widehat{\CC} \sm \ov{\DD}$
which satisfies $\varphi_c(\infty)=\infty$ and
$\varphi'_c(\infty)=1$. It conjugates the dynamics of $Q_c$ to
the squaring map so that $\varphi_c(z^2+c)=(\varphi_c(z))^2$ for
all $z$ in the basin of attraction of infinity.
The analytic curve $\rd(t):=\varphi^{-1}_c
\{ re^{2 \pi i t} : r>1 \}$ is called the {\it dynamic ray} at
angle $t \in \TT$. One immediately obtains $Q_c(\rd(t))=\rd(\db(t))$, where
$\db: t \mapsto 2t$ (mod 1) is the doubling map on $\TT$.
We say $\rd(t)$ {\it lands} at $z \in J_c$ if
$\lim _{r \to 1} \varphi^{-1}_c (re^{2 \pi i t})=z$.

By a {\it cycle} of $Q_c$ we simply mean
a periodic orbit $z \mapsto Q_c(z) \mapsto \cdots \mapsto
Q_c^{\circ n}(z)=z$.
The quantity $\lambda:=(Q_c^{\circ n})'(z)$ is
called the {\it multiplier} of this cycle. The cycle is {\it
attracting}, {\it repelling}, or {\it indifferent} if $|\lambda|<1$,
$|\lambda|>1$, or $|\lambda|=1$, respectively. An indifferent cycle is
{\it parabolic} if its multiplier is a root of unity. A quadratic
polynomial has always infinitely many cycles in the plane,
but at most one of them can be non-repelling.

A point $c \in \CC$ is called a {\it hyperbolic} parameter if
the sequence $\{ Q_c^{\circ n}(0) \}_{n \geq 0}$
tends to $\infty$ or to a necessarily unique attracting
cycle in $\CC$. It is called a {\it parabolic} parameter if $Q_c$ has a
necessarily unique parabolic cycle. Finally, $c$ is called a {\it Misiurewicz}
parameter if the critical point $0$ of $Q_c$ is preperiodic, i.e., if
$0$ has a finite forward orbit but is not periodic.

\subsection*{Hyperbolic components of $M$}
\label{subsec:hc}

We now turn to the parameter space.
Recall that the Mandelbrot set $M$ is the set of parameters $c$ for
which the (filled) Julia set of $Q_c$ is connected. Equivalently,
$c \in M$ if and only if $0 \in K_c$. Thus $Q_c^{\circ n}(0) \to \infty$ if
$c \notin M$, and it follows that all parameters outside $M$ are
hyperbolic. The hyperbolic parameters in $M$ form an open set;
in fact, each connected component of this set is a
connected component of the interior of $M$. As such, it is
called a {\it hyperbolic component} of $M$. The main
hyperbolic component $H_0$ containing $c=0$ is the prominently
visible cardioid in any picture of $M$. It consists of all $c$ for
which $Q_c$ has an attracting fixed point in $\CC$. The {\it period} of a
hyperbolic component $H$, denoted by $\per(H)$, is the length of
the unique attracting cycle of $Q_c$ for any $c \in H$. There is a
canonical conformal isomorphism $\lambda_H : H \iso \DD$ which
assigns to each $c \in H$ the multiplier of its attracting cycle.
The map $\lambda_H$ extends to a homeomorphism $\ov{H} \iso
\ov{\DD}$. The {\it center} and {\it root} of $H$ are by definition
the points $\lambda_H ^{-1}(0)$ and $\lambda_H ^{-1}(1)$,
respectively. The root of $H$ is a parabolic parameter; in fact every
parabolic parameter is realized as the root of a unique hyperbolic
component. There are exactly two parameter rays of angles
$\theta_-(H) < \theta_+(H)$ landing at the root of $H$ which are
rationals of the form $n/(2^p-1)$, where $p=\per(H)$. (When
$H=H_0$, the two angles $\theta_-(H_0)=0$ and $\theta_+(H_0)=1$
coincide, and only one ray $\rp(0)=\rp(1)$ lands at the root point
$c=1/4$.) It follows that these two angles have binary expansions
of the form
\begin{equation}
\label{eqn:angle}
\theta_-(H)= 0.\, \ov{\theta_0} \ \ \ \operatorname{and}\ \ \
\theta_+(H)= 0 . \, \ov{\theta_1},
\end{equation}
where $\theta_0$ and $\theta_1$ are binary words of length $p$, and
the bars indicate infinite repetition as usual. Conversely,
every parameter ray $\rp(t)$ for which $t$ is rational of odd
denominator lands at the root of a unique hyperbolic component.

Given any hyperbolic component $H$ and any
irreducible fraction $0< p/q <1$, there exists a
unique hyperbolic component $W$ which satisfies
$$\ov{H} \cap \ov{W} = \lambda_H^{-1}( e^{2 \pi i p/q} )
=\lambda_{W}^{-1}(1).$$
This $W$ is usually called the {\it $p/q$-satellite} of $H$.
%%%%%%%%%%%%%%%%%%%%%%%%%%%%%%%%%%%%%%%%
\begin{figure}[tp]
        \begin{center}
        \hspace*{-1in}
        \includegraphics[width=.6\textwidth,angle=-90]{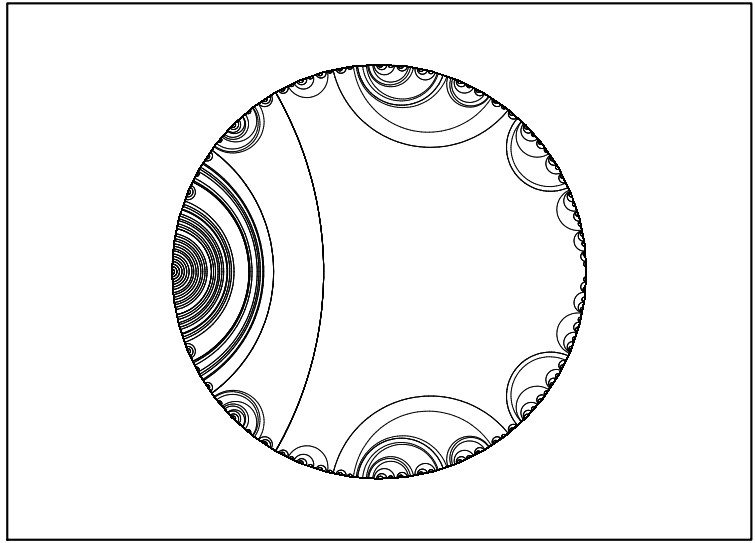}
        \end{center}
        \label{lam}
        \vspace{-1in}
        \caption{\sl The rational equivalence classes of
$\simeq$ used to define the abstract Mandelbrot set $M_{\operatorname{abs}}$.
The set $\cal R$ (see \S \ref{sec:quad}) corresponds to the closure of
the Poincar\'e geodesics which are symmetric with respect to the
real line, or equivalently, the ones which separate points $1$ and
$-1$ on the boundary circle.}
\end{figure}
%%%%%%%%%%%%%%%%%%%%%%%%%%%%%%%%%%%%%%%%

\subsection*{The abstract Mandelbrot set}
\label{subsec:Mabs}

Define an equivalence relation $\simeq$ on $\QQ / \ZZ$ by setting
$t \simeq s$ if and only if $\rp(t)$ and $\rp(s)$ land at a common
point. This can be extended to an equivalence relation on $\TT$ by
taking closure, and then to the closed disk $\ov{\DD}$ by taking
the Poincar\'e convex hulls of the equivalence classes on the
boundary. We still denote this equivalence relation by $\simeq$
(see Fig. 2). The {\it abstract Mandelbrot set} is by definition
the quotient $M_{\operatorname{abs}}:=\ov{\DD} / \! \! \simeq$.
It is a compact, full,
connected, and locally-connected space \cite{Douady2}.

The {\it MLC conjecture} asserts that the Mandelbrot set is
locally-connected, which, if true, would allow a complete topological
description of $M$. The celebrated {\it density of hyperbolicity conjecture}
asserts that all the interior components of $M$ are hyperbolic, so
that every quadratic polynomial can be approximated by a sequence
of hyperbolic quadratics. Douady and Hubbard have shown that MLC
implies density of hyperbolicity. In fact, they construct a
continuous surjection $\chi : M \to M_{\operatorname{abs}}$ whose fibers
$\chi^{-1}$(point) are reduced to points if and only if MLC holds.
In this case, the combinatorial model $M_{\operatorname{abs}}$ is actually
homeomorphic to $M$ via $\chi$. Density of hyperbolicity is equivalent to
the weaker statement that the fibers of $\chi$ have empty interior
(see \cite{Douady2}, \cite{Douady-Hubbard}, or \cite{Dierk}).

\subsection*{The tuning operation}
\label{subsec:tuning}

For every hyperbolic component $H$, denote by $\io_H : M
\hookrightarrow M$ the Douady-Hubbard's {\it tuning} map (see
\cite{Douady-Hubbard}, \cite{Haissinsky}, and
\cite{Milnor}). Then $\io_H$ is a topological embedding which maps
$H_0$ onto $H$ and respects centers and roots of hyperbolic
components. The image $\io_H(M)$ is what is often called the
{\it small copy of $M$ growing from $H$}. For a hyperbolic component $W$,
the image $\io_H(W)$ is a hyperbolic component of period $\per(H)
\cdot \per(W)$, which is called {\it $H$ tuned by $W$}. One has the
relation $\lambda_W(c)=\lambda_{\io_H(W)} (\io_H(c))$ for every $c
\in W$. The binary operation $(H,W) \mapsto \io_H(W)$ makes the set
of hyperbolic components into a free semigroup with $H_0$ as the
two sided identity.

The effect of this tuning map on Julia sets can be roughly
described as follows. Let $H$ be a hyperbolic component, $c_1 \in
M$ and $c_2:=\iota_H(c_1)$. Take the filled Julia set $K_{c}$ for
any $c \in H$ and replace each bounded component of $\CC \sm K_c$
by a copy of the filled Julia set $K_{c_1}$ by appropriately
identifying their Carath\'eodory loops. The resulting compact set
is homeomorphic to the filled Julia set $K_{c_2}$; see
\cite{Haissinsky} for details.

The tuning operation also acts on external angles. Let $\theta_-(H)
= 0. \, \ov{\theta_0} < \theta_+(H) = 0. \, \ov{\theta_1}$ be the
angles of the two parameter rays landing at the root of a
hyperbolic component $H \neq H_0$ as in (\ref{eqn:angle}). Take an
angle $t \in \TT$ with the binary expansion $0. \, t_1 t_2 t_3
\cdots$. Define the tuned angle
\begin{equation}
\label{eqn:tangle}
A_H(t) := 0 . \, \theta_{t_1} \theta_{t_2} \theta_{t_3} \cdots
\end{equation}
obtained by concatenating blocks of words of length $p=\per(H)>1$.
Note that under this {\it tuning algorithm} on angles, a dyadic
rational has two distinct images since it has two different
binary representations. It can be shown that if $c \in M$ is the
landing point of $\rp(t)$, then the tuned point $\io_H(c)$ is the
landing point of $\rp(A_H(t))$. The image
$A_H(\TT)$ is a self-similar Cantor set. In fact, for $i=0,1$ let
$\Lambda_i : \TT \to \TT$ be the map defined by
$\Lambda_i(0.t_1t_2t_3 \cdots):=0. \theta_i t_1t_2t_3 \cdots$. Then
$\Lambda_i$ is an affine contraction by a factor $2^{-p}$ and the
image $A_H(\TT)$ is precisely the invariant set generated by
$(\Lambda_0, \Lambda_1)$ \cite{Manning}. A standard computation then
shows that $A_H(\TT)$ has Hausdorff dimension $1/p < 1$ and
hence measure zero (see for example \cite{Mattila}).

\subsection*{Real hyperbolic components}
\label{subsec:real hc}

Let ${\cal H}$ denote the collection of all hyperbolic components
of $M$ which intersect the real line. Every $H \in {\cal H}$ is
invariant under the conjugation $c \mapsto \ov{c}$, and has its
center on the real line. If $H \cap \RR = ] c' ,
c [$, then $c=\lambda_H^{-1}(1)$ is the root of $H$ and
$c'=\lambda_H^{-1}(-1)$. It follows that $c$ is
the landing point of two parameter rays at angles
%%%%%%%%%%%%%%%%%%%%%%%%%%%%%%%%%%%%%%
\realfig{opening}{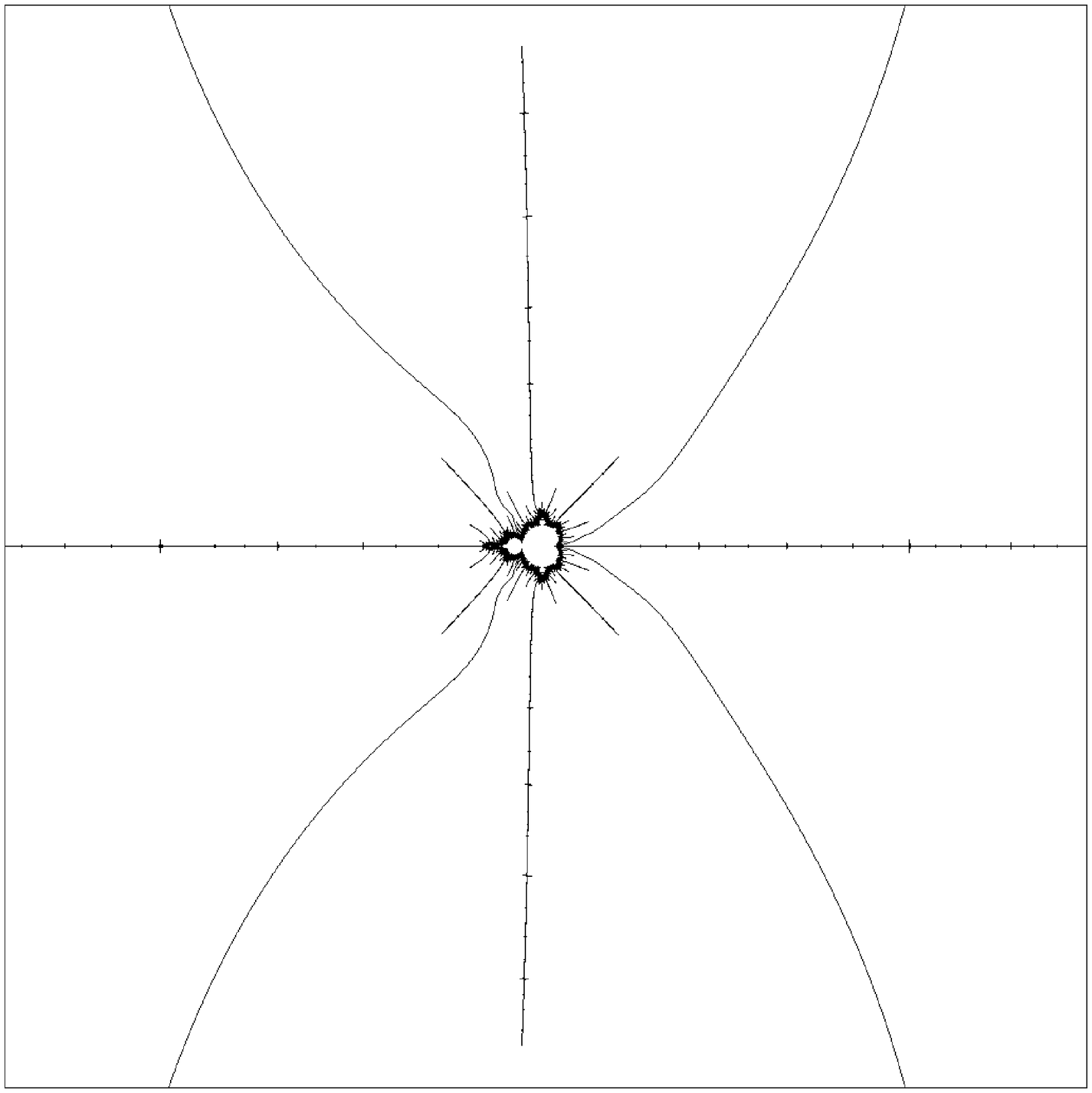}{\sl{A real period $5$ hyperbolic component
$H$ whose center is approximately located at $-1.985424253$. Here
$\theta_-(H)=15/31$ and $\omega_-(H)=16/33$. By definition, the
opening of this component is the interval $O(H)= \, ] 15/31 , 16/33
[$.}
$$
{
 \at{1.1\pscm}{9.0\pscm}{$\frac{15}{31}$}
 \at{-4.2\pscm}{9.0\pscm}{$\frac{16}{33}$}
 \at{1.1\pscm}{3.1\pscm}{$\frac{16}{31}$}
\at{-4.2\pscm}{3.1\pscm}{$\frac{17}{33}$}
}
$$}
{7cm}
%%%%%%%%%%%%%%%%%%%%%%%%%%%%%%%%%%%%%%
$$\theta_-(H)=\frac{n}{2^p-1} \ \ \ \ \operatorname{and}\ \ \ \
\theta_+(H)=1-\theta_-(H),$$
where $p=\per(H)$ and $n$ is an integer which satisfies
$0\leq n \leq 2^{p-1}-1$. A brief computation using the tuning algorithm shows
that $c'$ is the landing point of the parameter rays at angles
$$\omega_-(H)=\frac{n+1}{2^p+1} \ \ \ \ \operatorname{and}\ \ \ \
\omega_+(H)=1-\omega_-(H).$$
The open interval $O(H):=\, ]\theta_-(H), \omega_-(H)[ \, \subset
[0,1/2] $ is called the {\it opening} of $H$ (see \figref{opening}).
For example, $O(H_0)=\, ]0,1/3[$.

The following long-standing conjecture of Fatou, which has been
proved rather recently, will be used repeatedly in the next section
(see \cite{Lyubich1} and \cite{GS}):

\begin{theorem} [Density of real hyperbolics]
\label{real HD}
The set of all real hyperbolic parameters is dense in $M \cap \RR =
[-2,1/4]$. In particular, every component of the interior of $M$ which meets
the real line is hyperbolic.
\end{theorem}

%\vspace{1cm}

\section{Real quadratics and the set $\cal R$}
\label{sec:quad}

\subsection*{The $\tau$-function}
\label{subsec:tau}

Suppose $c \in \bd M \cap \RR$. By definition, the {\it dynamic
root} $r_c$ of $Q_c$ is the critical value $c$ if $c$ is not a
parabolic parameter so that the Julia set $J_c$ is full. On the other hand,
when $c$ is a parabolic parameter, the dynamic root $r_c$ is
the unique point of the parabolic cycle which is on the boundary
of the bounded Fatou component containing $c$.
In this case, $c<r_c$ and the open interval $] c , r_c [$
does not intersect $J_c$.

Recall that the {\it (prime-end) impression} of a parameter ray
$\rp(t)$ is the set of all $c \in \bd M$ for which there is a
sequence $\{ w_n \}$ such that $|w_n|>1$, $w_n \to e^{2 \pi i t}$, and
$\Phi^{-1}(w_n) \to c$. We denote the impression of
$\rp(t)$ by $\hrp(t)$. It is a non-empty, compact, connected subset
of $\bd M$. Every point of $\bd M$ belongs to the impression of at least
one parameter ray. Conjecturally, every parameter ray $\rp(t)$ lands at a
well-defined point $c(t) \in \bd M$ and $\hrp(t)= \{ c(t) \}$.
According to Douady-Hubbard and also Tan Lei, this certainly holds
for every rational angle $t$ and the landing point is either
parabolic or Misiurewicz depending on whether $t$ has
odd or even denominator (compare \cite{Douady-Hubbard},
\cite{Tan1} and \cite{Tan2}). Moreover, one can describe
which rational external rays land at parabolic and
Misiurewicz parameters. In the special case of real quadratics,
their result yields the following

\begin{theorem} [Douady-Hubbard-Tan\,Lei]
\label{dh}
Let $c \in \bd M \cap \RR$ be parabolic or Misiurewicz.
Then there exists a unique angle $\tau(c) \in [0,1/2]$ such that the dynamic
rays $\rd(\pm \tau(c))$ land at the dynamic root $r_c$ of $Q_c$. In the
parameter plane, the two rays $\rp(\pm \tau(c))$ land at $c$; in fact
$\hrp(\pm \tau(c))=\{ c \}$ and no other parameter
ray can have $c$ in its impression.
\end{theorem}

Note that the cases $c =1/4$ and $c=-2$ are special since
$\tau(1/4)=0$ and $\tau(-2)=1/2$ and the two rays given by the
theorem coincide. The following statement is immediate:

\begin{corollary}
\label{noreal}
If $t$ belongs to the opening $O(H)$ of some real hyperbolic component
$H$, then the impression $\hrp (t)$ of the parameter ray at angle $t$
does not intersect the real line.
\end{corollary}

It is not hard to partially generalize \thmref{dh} to all real quadratics
with connected Julia sets. The dynamic part of the following proof
which uses harmonic measure on the Julia set is inspired by the
more general combinatorial arguments in \cite{Zakeri}.

\begin{theorem}
\label{2rays}
Let $c \in \bd M \cap \RR$. Then there exists a unique angle
$\tau(c)\in [0,1/2]$ such that the dynamic rays $\rd( \pm \tau(c))$
land at the dynamic root $r_c$ of $Q_c$. In the parameter plane,
the two rays $\rp(\pm \tau(c))$, and only these rays, contain $c$
in their impression.
\end{theorem}

\figref{Feig} illustrates the content of this theorem.

\begin{pf}
In view of \thmref{dh} we may assume that $c$ is neither parabolic nor
Misiurewicz. First consider the dynamic plane. According to
\cite{Levin-VS}, the Julia set $J_c$ is locally-connected.
In particular, by the theorem of Carath\'eodory, all dynamic rays
land. By real symmetry, there exists at least one angle $0<t<1/2$
such that the two dynamic rays $\rd(\pm t)$ land at $r_c=c$. Assume
by way of contradiction that there are two angles $0<s<t<1/2$ such
that $\rd(s)$ and $\rd(t)$ both land at $c$. Consider the component
$W$ of $\CC \sm \ov{\rd(s) \cup \rd(t)}$ which does not intersect
the real line, and set $L_0:=\ov{W} \cap J_c$.
By the theorem of F. and M.~Riesz, $L_0$ is a non-degenerate continuum
of positive harmonic measure $\mu(L_0)=t-s$. Set $L_n:=Q_c^{\circ n}(L_0)$
and $c_n:=Q_c^{\circ n}(c)$.

%%%%%%%%%%%%%%%%%%%%%%%%%%%%%%%%%%%%%%%%%
\begin{figure}[tp]
  \hbox to \hsize{
   \hbox{\vbox{\hsize=.455\hsize \columnwidth=\hsize
          \centerline{\psfig{figure=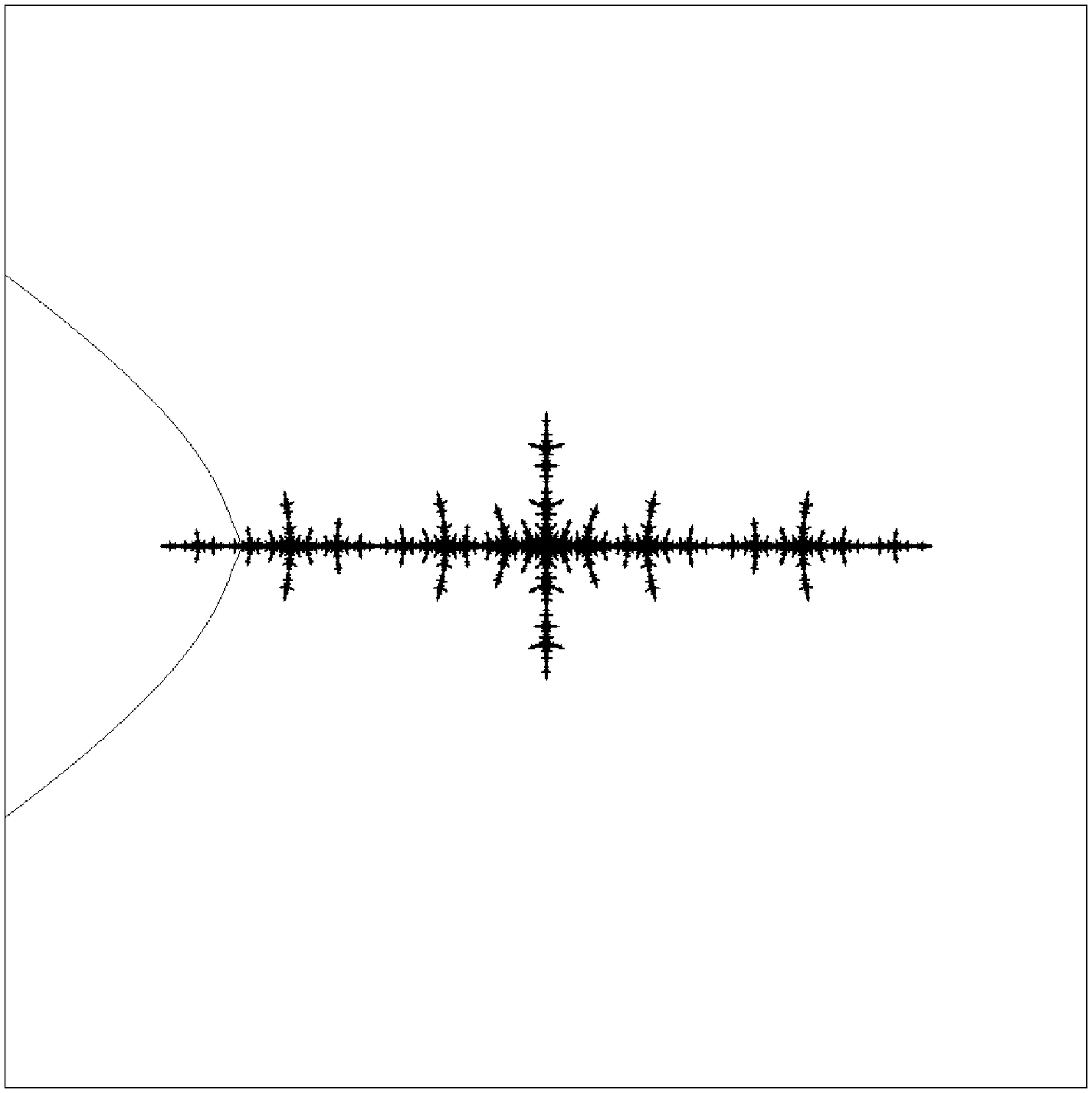,width=7cm}}
        }}
    \hfil
    \hbox{\vbox{\hsize=.455\hsize \columnwidth=\hsize
          \centerline{\psfig{figure=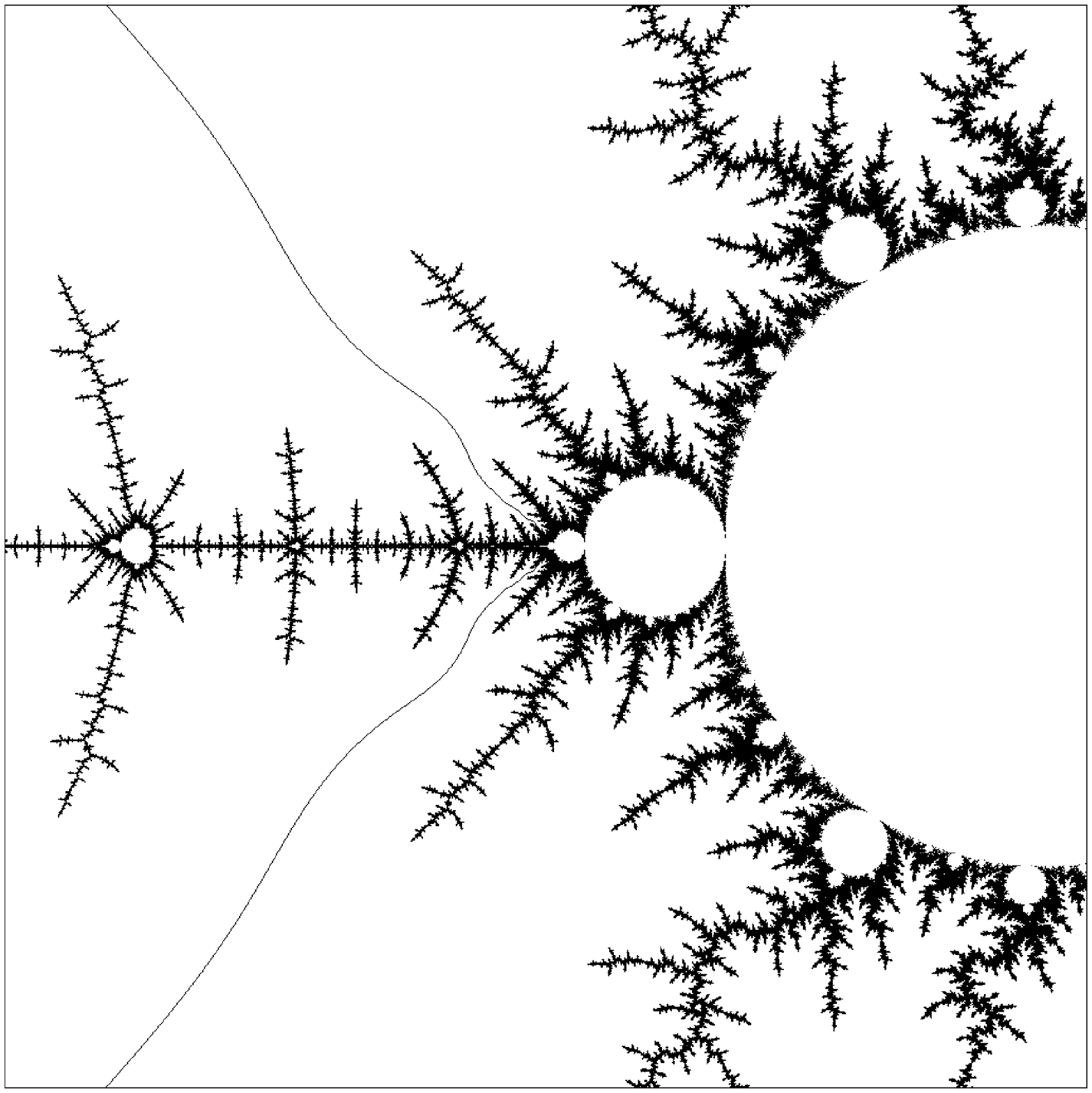,width=7cm}}
         }}
   }
\caption[]{{\sl Left: The Julia set of the main Feigenbaum
quadratic \mbox{$z \mapsto z^2 +c_{\operatorname{Feig}}$}, where
$c_{\operatorname{Feig}} \approx -1.401155$, with the dynamic
rays at angles $\pm \tau(c_{\operatorname{Feig}})$ landing at its
critical value. Computation gives $\tau(c_{\operatorname{Feig}})
\approx 0.412454$. Right: The corresponding parameter rays
outside the Mandelbrot set. }} \oplabel{Feig}
\end{figure}
%%%%%%%%%%%%%%%%%%%%%%%%%%%%%%%%%%%%%%%%%

It is easy to see that $J_c \cap \RR = [-\beta, \beta]$,
where the fixed point $\beta > 0$ is the landing point of $\rd(0)$.
Note that
$$Q_c^{-1}[-\beta, \beta]= [-\beta, \beta] \cup [i\xi, -i\xi],$$
where $\xi = \sqrt{ \beta+c } >0$ is the positive root of the
equation $Q_c(\pm i \xi) = - \beta$.  We claim that
\begin{equation} \label{eqn:LL}
L_n \cap \RR =\{ c_n \} \quad \text{for all} \ n \geq 0.
\end{equation}
Otherwise, take the smallest integer $n \geq 1$ for which this is false.
Then, since $L_{n-1} \cap \RR =\{ c_{n-1} \}$,
$L_{n-1}$ must intersect $[i\xi,-i\xi]$ at some point $z$. Since
$L_{n-1}$ is path-connected, it follows that it contains the unique
arc in $J_c$ which joins $z$ to $c_{n-1} \in \RR \sm \{ 0 \}$.
Evidently this arc has to intersect $\RR$ along a
non-degenerate interval, contradicting $L_{n-1} \cap \RR =\{ c_{n-1} \}$.
Thus \eqref{eqn:LL} holds.

Now it easily follows that each map $Q_c: L_n \mapsto L_{n+1}$ is
injective. In fact, if this were not true for some $n\geq 0$,
then $L_n$ would have to contain a pair of symmetric
points $\pm z$ in $J_c$ and hence
the critical point $0$. Since $L_n \cap \RR =\{ c_n \}$, this would mean
$c_n=0$, which would be impossible. Thus, in the cascade of injective maps
$L_0 \mapsto L_1 \mapsto \cdots \mapsto L_n \mapsto \cdots $,
the harmonic measure $\mu$ doubles at each step. However,
$\mu(L_0)>0$ and $\mu(L_n) \leq \mu (J_c)=1$ for all $n$. The contradiction
shows there is a unique $t=\tau(c) \in [0,1/2]$ such that $\rd(t)$
lands at $c$, and the proof in the dynamic plane is complete.

To prove the result in the parameter plane, one possible approach is
to use the following standard construction:
For any $c \in \bd M \cap \RR$, the
{\it itinerary} of the dynamic root $r_c$ is the infinite sequence
${\bf itin}(r_c):=(\ve_0, \ve_1,\ve_2 , \ldots)$ of signs
$\pm$ determined by $\ve_j=\operatorname{sgn}(Q_c^{\circ j}(r_c))$.
It is easy to see that the binary expansion $0.\,t_0\, t_1\, t_2\cdots$
of the angle $\tau (c)$ constructed above is uniquely determined by
${\bf itin}(r_c)$. In fact, $t_0=0$ and a brief computation shows that
for $j \geq 0$,
$$t_{j+1}= \left \{ \begin{array}{cl}
t_j & \quad \text{if} \quad\ve_j = + \vs \\
1-t_j & \quad \text{if} \quad \ve_j = -
\end{array} \right. $$
It follows in particular that
if ${\bf itin}(r_c)$ and ${\bf itin}(r_{c'})$ coincide up to
the first $n$ signs, then $|\tau(c)-\tau(c')| \leq 2^{-(n+1)}$.

Now consider a parameter $c \in \bd M \cap \RR$ which is neither
parabolic nor Misiurewicz, and let $0<s<1/2$ be any angle for which
$c \in \hrp(s)$. Choose a decreasing (resp. increasing) sequence
$\{ a_n \}$ (resp. $\{ b_n \} )$ of real parabolic parameters
converging to $c$ (the existence of such sequences is guaranteed
by \thmref{real HD}). It follows from \thmref{dh} that
\begin{equation}
\label{eqn:ss}
\tau (a_n) < \tau (a_{n+1}) < s < \tau (b_{n+1}) < \tau (b_n)
\end{equation}
for all $n$. It is not hard to see that the
sequences ${\bf itin}(r_{a_n})$ and ${\bf itin}(r_{b_n})$ converge to
${\bf itin}(c)$ in the 2-adic metric. Hence both $\tau(a_n)$ and
$\tau (b_n)$ converge to $\tau(c)$ and it
follows from (\ref{eqn:ss}) that $s=\tau(c)$.
Since $c$ must belong
to the impression of at least one parameter ray, it follows that $c
\in \hrp (\tau(c))$ and this completes the proof.
\end{pf}

\subsection*{The set $\cal R$}
\label{subsec:R}

We now consider the set of external angles in
$[0,1/2]$ of the parameter rays whose impressions intersect the real
line:
\begin{equation}
\label{eqn:R}
{\cal R}:=\{ t \in [0,1/2] : \hrp(t) \cap \RR \neq \es \}.
\end{equation}

\begin{lemma}
\label{R=c} The mapping $\tau : \bd M \cap \RR \to {\cal R}$ given
by \thmref{2rays} is a homeomorphism. Its inverse $\pi=\tau^{-1}$
is determined by $\{ \pi(t) \} = \hrp(t) \cap \RR$.
\end{lemma}

\begin{pf}
Given $t \in {\cal R}$, let us show that the impression
$\hrp(t)$ intersects the real line at a unique point. By
\thmref{dh} this is true if $t=0$ or $1/2$, so we may assume
\mbox{$0< t <1/2$}. Let $c,c' \in \bd M \cap \RR$ both belong to
$\hrp(t)$, and $c'<c$. By symmetry, $c$ and $c'$ belong to the impression
$\hrp(-t)$ also. By \thmref{real HD}, there exists a real hyperbolic component
$H$ such that $H \cap \RR \subset \, ]c',c[$. It follows that
the union $\hrp(t) \cup \hrp(-t)$ separates the plane and $H$ is contained in
a bounded component $W$ of $\CC \sm (\hrp(t)\cup \hrp(-t))$.
Now $\bd W \subset \hrp(t)\cup \hrp(-t) \subset \bd M$, so
$W$ must be a component of the interior of $M$, implying $W=H$.
In particular, the root of $H$ belongs to $\hrp(t)$.
This, by \thmref{dh}, implies $\hrp(t)$ is a singleton, which
contradicts our assumption.

Thus the map $\pi:{\cal R} \to \bd M \cap \RR$ given by
$\hrp(t) \cap \RR=\{ \pi(t) \}$ is well-defined. The relations $\pi
\circ \tau =\operatorname{id}$ and $\tau \circ
\pi=\operatorname{id}$ follow easily from \thmref{2rays}. In
particular, $\tau$ is both injective and surjective.

It remains to prove continuity of $\tau$, or equivalently $\pi$.
Clearly $\pi$ is monotone. Assume by way of contradiction that $\{
t_n \}$ is an increasing sequence in $\cal R$ converging to $t \in
{\cal R}$ such that $c^{\ast}:=\lim \pi(t_n) > c:=\pi(t)$. If there
exists some $c' \in \bd M \cap \RR$ in the interval $]c,c^{\ast}[$,
then $\tau(\pi(t_n)) < \tau(c') < \tau(c)$ or $t_n < \tau(c') < t$
for all $n$, which is impossible. Hence $]c,c^{\ast}[$ is a subset
of the interior of the Mandelbrot set. By the density of real
hyperbolics and the fact that $c,c^{\ast} \in \bd M$,
there exists a component $H \in {\cal H}$ such that $H \cap
\RR =]c,c^{\ast}[$. Take, for example, the root point
of the $1/3$-satellite of $H$ which is the landing point of two rational
parameter rays at angles $\alpha < \beta$. Then $t_n < \alpha <
\beta < t$ for all $n$, which again is a contradiction. This proves
that $\pi$ is left-continuous. The proof of right-continuity is
similar.
\end{pf}

\begin{lemma}
\label{Odense}
The union $\bigcup_{H \in {\cal H}} O(H)$ of the
openings of real hyperbolic components is dense in $[0,1/2]$.
\end{lemma}

\begin{pf}
Let $E$ be the set of endpoints of the openings $O(H)$ for
$H \in {\cal H}$. Evidently, $E \subset {\cal R}$. Since
$\cal R$ is closed by \lemref{R=c}, we actually have $\ov{E} \subset
{\cal R}$.

Assuming the lemma is false, let $[t,t']$ be a maximal non-degenerate
interval in the complement of the union $\bigcup_{H \in {\cal H}}
O(H)$. Then $t,t' \in \ov{E}$, so both $t$ and $t'$ are angles in
$\cal R$. Applying the homeomorphism $\pi$, we obtain $c':=\pi(t') <
c:=\pi(t)$. Density of real hyperbolics now implies the existence
of some $H \in \cal H$ with $H \cap \RR \subset \, ]c',c[$. It follows that
$]t,t'[$ must contain the opening $O(H)$, which is a
contradiction.
\end{pf}

\begin{corollary}
\label{R=op}
${\cal R}=[0,1/2] \sm \bigcup_{H \in {\cal H}} O(H)$.
\end{corollary}

\begin{pf}
The inclusion $\subset$ follows from \corref{noreal}. To see the
inclusion $\supset$, note that by \lemref{Odense}, any $t \in [0,1/2]$
outside the union $\bigcup_{H \in {\cal H}} O(H)$ belongs to $\ov{E}$,
which is a subset of ${\cal R}$.
\end{pf}

The definition of ${\cal R}$ is simple but rather hard to work
with. In fact, by \corref{R=op} the explicit construction of ${\cal
R}$ boils down to deciding which rational angles form the
endpoints of the openings of real hyperbolic components.
There is a combinatorial algorithm due
to P.~Lavaurs which describes the equivalence relation $\simeq$ on
$\QQ/ \ZZ$ used in the definition of $M_{\operatorname{abs}}$ in \S
\ref{sec:bg} \cite{Lavaurs}. In other words, it tells which
rational angles $t,s$ satisfy $t \simeq s$, or equivalently, which
rational ray pairs $(\rp(t),\rp(s))$ land at a common point. The
rational parameter rays landing on the real line correspond to the
angles $t$ in $[0,1/2]$ for which $t \simeq -t$. The set ${\cal R}$ is obtained
by taking the closure of the set of such rational angles. Thus,
{\it in principle}, one should be able to determine the set $\cal
R$ using Lavaurs' algorithm. However, this algorithm is not quite
suitable for the purpose of computing measure and dimension of
$\cal R$. To circumvent this problem, we found an alternative
description of $\cal R$, much easier to work with since it is
purely given by the dynamics of the doubling map \mbox{$\db: t
\mapsto 2t$ (mod 1)} on the circle. As it turned out, this
description was not new: In \cite{Douady3}, Douady gives
a description for the dynamic rays which land on the spine
$[-\beta, \beta]$ of the Julia set of a real quadratic polynomial
(compare equation (\ref{eqn:S_c}) below). When transferred to the
parameter plane, it would give the same alternative description of
$\cal R$.

\begin{theorem}
\label{R=E}
${\cal R} = \{ t \in [0,1/2] : \db ^{\circ n}(t) \notin \, ] t,1-t [ \
\operatorname{for\ all\ } n\geq 1 \} $.
\end{theorem}

\begin{pf}
Suppose that $t \in {\cal R}$ and $c=\pi(t)$, or equivalently
$t=\tau(c)$ as in \thmref{2rays}. This means the dynamic rays $\rd(\pm t)$ land
at the dynamic root $r_c$ of $Q_c$, and hence the dynamic rays
$\rd(\db^{\circ n}(\pm t))$ land at $Q_c^{\circ n}(r_c)$ for all $n\geq 1$.
Note that since $c$ and $r_c$ are real, $Q_c^{\circ
n}(r_c) \geq c$ for all $n \geq 1$. Moreover, since
$J_c \, \cap \, ] c, r_c [ = \emptyset$,
we actually have $Q_c^{\circ n}(r_c) \geq
r_c$. It easily follows that $\db^{\circ n}(t) \notin \, ] t,1-t [$
for all $n \geq 1$.

Now let $t \notin {\cal R}$. By \corref{R=op}, this means $t \in
O(H)$ for a real hyperbolic component $H$ of some period $p \geq 1$,
so that
$$\frac{n}{2^p-1}=\theta_-(H) < t < \omega_-(H)=\frac{n+1}{2^p+1}.$$
Thus $t<2^pt-n<1-t$, or $\db^{\circ p}(t) \in \, ]t,1-t[$.
\end{pf}

\section{Measure and dimension of $\cal R$}
\label{sec:AB}

The alternative description of $\cal R$ given by \thmref{R=E}
makes the question of measure of $\cal R$ almost trivial:

\begin{lemma}\label{ergodic}
$\cal R$ is a set of Lebesgue measure zero.
\end{lemma}

\begin{pf}
Choose a nested sequence $I_1 \supset I_2 \supset \cdots$ of open
intervals centered at $1/2$ such that $\bigcap_{n \geq 1} I_n = \{
1/2 \}$. Since $\db$ is ergodic with respect to Lebesgue measure,
for each $n$ there exists
a set $X_n \subset \TT$ with $m(X_n)=0$ such that the forward orbit
of $t$ under $\db$ hits $I_n$ whenever $t \in \TT \sm X_n$. Taking
$X:= \bigcup _{n \geq 1} X_n$, it follows that the forward orbit of
every $t \in \TT \sm X$ hits every $I_n$. Since for every $t \in
[0,1/2[$ there exists a large $n$ with $I_n \subset \, ] t , 1-t
[$, it follows that the forward orbit of every $t \in \TT \sm X$
hits $]t,1-t[$. By \thmref{R=E}, we must have ${\cal R} \subset X$,
which proves the claim.
\end{pf}

\noindent
{\it Proof of \thmref{A}.} The set of external angles of
the parameter rays which land on the real slice is evidently
contained in ${\cal R} \cup -{\cal R}$, and thus has measure
zero by \lemref{ergodic}. It follows that $\mu_M [-2,1/4]=0$.
$\hfill \Box$ \vspace{2mm}

Using \corref{R=op}, we obtain

\begin{corollary}\label{sum}
The sum of the lengths of the openings of all real
hyperbolic components is $1/2$.
\end{corollary}

We now turn to the question of dimension of $\cal R$. First let us
introduce some notation. By a {\it dyadic rational of generation}
$n$ in $\TT = \RR / \ZZ$ we mean a number of the form $p/2^n$
modulo $\ZZ$, where $p$ is an odd integer. The set of all such
numbers will be denoted by ${\cal D}_n$. For $t\in \TT$ and $n\geq
1$, the notation $\| t \| _n$ will be used for the distance from
$t$ to the closest dyadic rational of generation $n$:
$$\| t \| _n := \inf_{x \in {\cal D}_n} |t-x|.$$
Clearly $0\leq \| t \|_n \leq 2^{-n}$ so that $\|
t \|_n \to 0$ as $n \to \infty$. For $0 < \sigma < 1$ and $n \geq
2$, we consider the nested sequence of non-empty compact sets
$${\cal K}^n_{\sigma}:= \{ t \in \TT : \| t \|_k \geq \sigma \,
2^{-k} \ \operatorname{for\ all\ } 2\leq k \leq n \},$$
and we define
\begin{equation}
\label{eqn:EE}
{\cal K}_{\sigma}:=\bigcap_{n\geq 2} {\cal K}^n_{\sigma} =
\{ t \in \TT : \| t \|_n \geq \sigma \,
2^{-n} \ \operatorname{for\ all\ } n \geq 2 \}.
\end{equation}

The proof of \thmref{B} will depend on the following

\begin{lemma}
\label{E_sigma}
$$\lim_{\sigma \to 0}\ \Hdim ({\cal K}_{\sigma}) = \lim_{\sigma \to 0}
\ \Hdim \left( {\cal K}_{\sigma} \cap \left[
\frac{1-\sigma}{2},\frac{1}{2} \right] \right) = 1.$$
\end{lemma}

Before proving this lemma, let us show how \thmref{B} would follow.
\vspace{2mm}\\
{\it Proof of \thmref{B} (assuming
\lemref{E_sigma}).} First observe that for all small $\sigma >0$,
\begin{equation} \label{eqn:rk}
{\cal R} \supset {\cal K}_{\sigma} \cap \left[ \frac{1-\sigma}{2} ,
\frac{1}{2} \right] .
\end{equation}
In fact, if $t \in {\cal K}_{\sigma}$,
then $\| t \| _{n+1} \geq \sigma 2^{-(n+1)}$ for all $n \geq 1$.
Applying the iterate $\db^{\circ n}$, it follows that $|\db^{\circ
n}(t) - 1/2| \geq \sigma /2$ for all $n \geq 1$. If in addition $t \in [
(1-\sigma)/2 , 1/2 ]$, then $\sigma /2 \geq 1/2 -t$, so that
$|\db^{\circ n}(t) - 1/2| \geq 1/2 -t$ for all $n \geq 1$. By
\thmref{R=E}, this implies $t \in {\cal R}$, which proves
(\ref{eqn:rk}). Thus,
$$1\geq \Hdim ({\cal R}) \geq \Hdim \left( {\cal K}_{\sigma} \cap
\left[ \frac{1-\sigma}{2} , \frac{1}{2} \right] \right) .$$
Taking the limit as $\sigma \to 0$, we obtain $\Hdim ({\cal R})=1$.
But by the theorem of \mbox{Beurling}, after removing a set of zero
capacity (and hence zero Hausdorff dimension) from $\cal R$, we may
assume that all the remaining rays in $\cal R$ land at a real
parameter. This completes the proof of \thmref{B}. \hfill $\Box$
\vspace{2mm}

The proof of \lemref{E_sigma} will be based on the following two
lemmas. The first one describes the metric structure of the compact sets
${\cal K}^n_{\sigma}$, and the second one uses this structure to
estimate the dimension of ${\cal K}_{\sigma}$.
For simplicity we state these lemmas in the
case $\sigma$ is a negative power of $2$.

\begin{lemma} \label{structure}
Fix a parameter $\sigma=2^{-p}$ where $p \geq 2$ is an integer, and
let $n\geq 2$. Then
\begin{enumerate}
\item[(i)]
${\cal K}^n_{\sigma}$ is the disjoint union of a finite collection $|{\cal
K}^n_{\sigma}|$ of closed non-degenerate intervals in $\TT$
whose endpoints are dyadic rationals in $\bigcup_{k=2}^n {\cal
D}_{k+p}$.
\item[(ii)]
There are two distinguished intervals $I^n_0 , I^n_{1/2} \in
|{\cal K}^n_{\sigma}|$ centered at $0$ and $1/2$, with
$m(I^n_0)=m(I^n_{1/2})=2^{-n+1}(1-\sigma)$.
\item[(iii)]
Every $I \in |{\cal K}^n_{\sigma}|$ contains at least one and at
most three elements of $|{\cal K}^{n+1}_{\sigma}|$. If $I=I^n_0$ or
$I^n_{1/2}$, then it contains exactly 3 elements of $|{\cal
K}^{n+1}_{\sigma}|$.
\item[(iv)]
For every $I \in |{\cal K}^n_{\sigma}|$,
$0 < m(I) \leq 2^{-n+1}(1-\sigma)$.
\item[(v)]
For every $I \in |{\cal K}^n_{\sigma}|$,
$$\frac{m(I \cap {\cal K}^{n+1}_{\sigma})}{m(I)} \geq
\frac{3-8\sigma}{3-4\sigma}.$$
\end{enumerate}
\end{lemma}

\begin{pf}
For convenience, let us introduce the sets
$${\cal E}^k_{\sigma}:= \{ t \in \TT : \| t \| _k \geq \sigma 2^{-k}
\}, \ \ \ \ \ k=2,3, \ldots $$
Evidently, each ${\cal E}^k_{\sigma}$ is the union of a collection
$|{\cal E}^k_{\sigma}|$ of $2^{k-1}$ disjoint closed intervals of
length $2^{-k+1}(1-\sigma)$ whose endpoints belong
to ${\cal D}_{k+p}$.

To prove (i), note that
$${\cal K}^n_{\sigma} = {\cal E}^2_{\sigma} \cap {\cal E}^3_{\sigma}
\cap \cdots \cap {\cal E}^n_{\sigma},$$ so the statement about the
endpoints of the intervals in $|{\cal K}^k_{\sigma}|$ follows from
the corresponding statement for $|{\cal E}^k_{\sigma}|$. If some $I
\in |{\cal K}^n_{\sigma}|$ degenerate to a singleton $\{ x \}$,
then $x$ must be the common endpoint of two intervals in $|{\cal
E}^k_{\sigma}|$ and $|{\cal E}^s_{\sigma}|$ for some $2\leq k < s
\leq n$. This is clearly impossible since ${\cal D}_{k+p} \cap
{\cal D}_{s+p} = \emptyset$.

The statements (ii) and (iii) are trivial. The statement (iv) follows
from the fact that every interval $I \in |{\cal K}^n_{\sigma}|$ is
contained in an interval $J \in |{\cal E}^n_{\sigma}|$, so that
$m(I) \leq m(J) \leq 2^{-n+1}(1-\sigma)$.

It remains to prove (v). If $I=I^n_0$ or $I^n_{1/2}$, a brief
computation shows that
$$\frac{m(I \cap {\cal K}^{n+1}_{\sigma})}{m(I)}=
\frac{1-2\sigma}{1-\sigma} \geq \frac{3-8\sigma}{3-4\sigma}.$$
So take an $I \in |{\cal K}^n_{\sigma}|$ such that $I \neq I^n_0$ and
$I \neq I^n_{1/2}$. Let
$$[ x+\sigma 2^{-n}, x+2^{-n+1}-\sigma 2^{-n}] =
[ x+2^{-(n+p)}, x+2^{-n+1}-2^{-(n+p)}]$$
be the unique interval in $|{\cal E}^n_{\sigma}|$ which contains $I$, where
$x \in {\cal D}_n$. Let
$$t:= x+2^{-(n+1)}, \ y:=x+2^{-n}.$$
Clearly $t \in {\cal D}_{n+1}$ and $y \in {\cal D}_m$ for some $m <
n$. It is easy to see that one of the following two cases must
occur:\vs

{\it Case 1.} $m+p=n+1$ and $I = [ x+2^{-(n+p)} , y-2^{-(n+1)}] =[
x+2^{-(n+p)} , t ]$. In this case
$$\frac{m(I \cap {\cal K}^{n+1}_{\sigma})}{m(I)}=
1-\frac{2^{-(n+p+1)}}{2^{-(n+1)}-2^{-(n+p)}}=\frac{1-3\sigma}{1-2\sigma}
\geq \frac{3-8\sigma}{3-4\sigma}.$$

{\it Case 2.} $m+p > n+1$ and $I = [ x+2^{-(n+p)} , y-2^{-(m+p)}]$.
In this case, $t$ belongs to the interior of $I$ and we have
$$\frac{m(I \cap {\cal K}^{n+1}_{\sigma})}{m(I)}=
1-\frac{2^{-(n+p)}}{2^{-n}-2^{-(n+p)}-2^{-(m+p)}} \geq
1-\frac{2^{-(n+p)}}{2^{-n}-2^{-(n+p)}-2^{-(n+2)}}
= \frac{3-8\sigma}{3-4\sigma}.$$
In either case, we obtain the lower bound in (v).
\end{pf}

\begin{lemma}
\label{aux}
Fix a parameter $\sigma=2^{-p}$ where $p \geq 2$ is an integer, and let
$$ 1 < \lambda:=\frac{3-4\sigma}{3-8\sigma} \leq 2.$$
Then
$$\Hdim({\cal K}_{\sigma}) \geq \Hdim \left( {\cal K}_{\sigma} \cap \left[
\frac{1-\sigma}{2},\frac{1+\sigma}{2} \right] \right) \geq
1-\frac{\log \lambda}{\log 2}. \vs$$
\end{lemma}

Note that \lemref{E_sigma} follows immediately, since by symmetry
$$\Hdim \left( {\cal K}_{\sigma} \cap \left[
\frac{1-\sigma}{2},\frac{1+\sigma}{2} \right] \right) =
\Hdim \left( {\cal K}_{\sigma} \cap \left[
\frac{1-\sigma}{2},\frac{1}{2} \right] \right)$$
and $\lambda \to 1$ as $\sigma = 2^{-p} \to 0$.

\begin{pf}
This follows from \lemref{structure}(v) by a standard mass distribution argument.
Define for each $n$ the probability measure $\mu_n$
supported in ${\cal K}^n_{\sigma}$, with uniform density on each interval in
$|{\cal K}^n_{\sigma}|$, as follows:
For $I \in |{\cal K}^2_{\sigma}|$, set $\mu_2(I):=m(I)/
m({\cal K}^2_{\sigma})$. When $n \geq 2$ and $I \in
|{\cal K}^{n+1}_{\sigma}|$, let $J$ be the unique interval in
$|{\cal K}^{n}_{\sigma}|$ that contains $I$, and set
$$\mu_{n+1}(I):=\frac{m(I)}{m(J \cap {\cal K}^{n+1}_{\sigma})} \,
\mu_n(J).$$
By \lemref{structure}(v),
$$\frac{\mu_{n+1}(I)}{m(I)}=\frac{m(J)}{m(J \cap {\cal
K}^{n+1}_{\sigma})} \cdot \frac{\mu_n(J)}{m(J)} \leq \lambda
\, \frac{\mu_n(J)}{m(J)}.$$
Continuing inductively, it follows that
$$\frac{\mu_{n+1}(I)}{m(I)} \leq (\con)\, \lambda ^{n}.$$
Now let $\mu^\ast$ be the weak limit of the sequence $\{ \mu_n
\}$. Then $\mu^\ast$ is supported in ${\cal K}_{\sigma}$ and we have
\begin{equation}
\label{eqn:lam}
\mu^\ast (I) \leq (\con)\, \lambda^{n}\, m(I)\quad \text{for all} \quad
I \in |{\cal K}^{n}_{\sigma}|.
\end{equation}
To estimate the $\mu^\ast$-measure of an arbitrary interval $T$ of
length $\ve = m(T) > 0$, choose $n$ so that $2^{-(n+1)} < \ve
\leq 2^{-n}$, and consider the union ${\cal I}$ of all the
intervals in $|{\cal K}^n_{\sigma}|$ which intersect $T$. Then,
by (\ref{eqn:lam}) and \lemref{structure}(iv),
$$\begin{array}{rl}
\mu^\ast(T) \leq \mu^\ast ({\cal I}) \leq & (\con) \, \lambda^{n}
\, m({\cal I}) \vspace{2mm} \\
 \leq & (\con)\, \lambda^{n} (\ve+2 \cdot 2^{-n+1})\vspace{2mm} \\
 \leq & (\con)\, \lambda^{n}\, \ve \vspace{2mm} \\
 \leq & (\con)\, \ve^{s} (\lambda^{n} 2^{-n(1-s)})
\end{array}$$
for all $s>0$. If $0<s<1-(\log \lambda / \log 2)$, we have
$\lambda  2^{s-1} < 1$ and hence
\begin{equation}\label{eqn:FR}
\mu^\ast(T) \leq (\con)\, m(T)^s.
\end{equation}
To finish the argument, let us recall the following (see for example
\cite{Mattila}): \vs
\begin{quotation}
{\bf Frostman's Lemma.} {\it A Borel set $X \subset \RR^d$ satisfies
$\Hdim(X)\geq s$ if and only if there exists a finite
Borel measure $\mu$ supported in $X$ and a constant $C>0$ such that
$\mu(B(x,r)) \leq C r^s$ for all $x \in \RR^d$ and all $r>0$.}
\end{quotation}
Applying this lemma to (\ref{eqn:FR}), we obtain the inequality
\begin{equation}
\label{eqn:lb}
\Hdim ({\cal K}_{\sigma}) \geq 1 - \frac{\log \lambda}{\log 2}.
\end{equation}

The claim for $\Hdim ( {\cal K}_{\sigma} \cap [ \,
 (1-\sigma)/2 , (1+\sigma)/2  \,  ] )$ follows from the same
argument: Just start with the natural probability measure on
${\cal  K}^{p+1}_{\sigma} \cap [ \,  (1-\sigma)/2 , (1+\sigma)/2 \, ]$,
and inductively define the sequence of measures $\mu_n$ supported in
${\cal K}^{n}_{\sigma} \cap [ \,  (1-\sigma)/2 , (1+\sigma)/2 \, ]$
for all $n > p+1$.
\end{pf}

The following observation will be used in the proof of \thmref{D}
in \S \ref{sec:bi}. Recall that a compact set $K$ of real numbers
is {\it porous} if there exists an $0< \ve < 1$ such that every
interval $I$ has a subinterval $J$ disjoint from $K$, with $m(J) >
\ve \, m(I)$. It is an easy exercise to show that a porous set has
Hausdorff dimension less than $1$.

\begin{lemma}
\label{<1}
For every $0< \sigma <1$, the compact set ${\cal K}_{\sigma}$ of
(\ref{eqn:EE}) is porous, hence $\Hdim ( {\cal K}_{\sigma} ) <1$.
\end{lemma}

\begin{pf}
Given an interval $I$ of length $m(I)>0$, choose $n$ so that $2^{-(n-1)} <
m(I) \leq 2^{-(n-2)}$ and consider a subinterval $] a/2^n ,
(a+1)/ 2^n [ \,\subset I$, with $a \in \ZZ$. Then the midpoint
$(2a+1)/2^{n+1}$ is dyadic of generation $n+1$, so the open interval
$$J:=\, \left ] \frac{2a+1-\sigma}{2^{n+1}} , \frac{2a+1+\sigma}{2^{n+1}}
\right [ \, \subset I$$
is disjoint from ${\cal K}_{\sigma}^{n+1}$ hence from
${\cal K}_{\sigma}$. Moreover,
$$
m(J) = \sigma \, 2^{-n} = \frac{\sigma}{4} 2^{-(n-2)} \geq
\frac{\sigma}{4} m(I).
$$
This proves ${\cal K}_{\sigma}$ is porous.
\end{pf}

\begin{remark}
When $\sigma=2^{-p}$, one can interpret $\cal K_{\sigma}$
symbolically as the set of all binary angles $0.t_1t_2t_3\cdots$ in $\TT$
such that if $t_jt_{j+1}\cdots t_{j+p-1}$ is the first occurrence of $p$
consecutive $0$'s or $1$'s, then $j=1$ or $j=2$. This description allows
a more combinatorial approach to the fact that
the Hausdorff dimension of $\cal K_{\sigma}$ tends to $1$ as
$\sigma$ tends to $0$.
\end{remark}

The rest of this section will be devoted to the proof of
\thmref{C}, which generalizes \thmref{B} to all tuned images of the
real slice. Consider a hyperbolic component $H \neq H_0$ with the
associated tuning maps $\io _H : M \hookrightarrow M$ on the
parameter plane and $A_H: \TT \rightarrow \TT$ on the external
angles (see \S \ref{sec:bg} and recall that $A_H$ is 2-valued at
every dyadic rational). Let $\eta_H:=\io_H [-2, 1/4]$ be the tuned
image of the real slice of $M$. Then, if $\rp(t)$ lands at $c \in
[-2,1/4]$, $\rp (A_H(t))$ lands at $\io_H(c) \in \eta_H$. Of
course, not every parameter ray landing on $\eta_H$ comes from
tuning. Along $\eta_H$ there are many ``branch points'' at which
more than two parameter rays land. All such branch points are
Misiurewicz parameters and hence the external angles of the rays
landing on them are rational with even denominator (a countable
set). It follows that such rays can be completely ignored when it
comes to computing Hausdorff dimensions.

Recall from \S \ref{sec:bg} that the image $K=A_H(\TT)$ is a
Cantor set of Hausdorff dimension $1/\per(H)$ which is invariant under the two
contractions $\Lambda_0$ and $\Lambda_1$. It is easy to see that
$A_H : [0,1] \to [\theta_-(H), \theta_+(H)]$ is
the right inverse of a ``devil's staircase.'' More precisely, there exists a
continuous non-decreasing map $\psi_H: [\theta_-(H), \theta_+(H)] \to  [0,1]$
which maps every gap of $K$ to a well-defined dyadic rational, mapping
$K$ onto $[0,1]$, such that $\psi_H \circ A_H(t)=t$ for all $t$.

The proof of \thmref{C} will depend on the following

\begin{lemma}
\label{Holder}
The map $\psi_H : [\theta_-(H), \theta_+(H)] \to  [0,1]$ is
H\"{o}lder continuous of exponent $1/\per(H)$.
\end{lemma}

\begin{pf}
Let $p=\per(H)$, $\theta_-(H)=0.\ov{\theta_0}$, and
$\theta_+(H)=0.\ov{\theta_1}$ as in (\ref{eqn:angle}). Pick any two points
$a<b$ in $[\theta_-(H), \theta_+(H)]$. We may assume that
$\psi_H(a) < \psi_H(b)$ since otherwise $\psi_H(a)=\psi_H(b)$ and
there is nothing to prove. Let $a' := \inf (K \cap [a, \theta_+(H)])$
and similarly define $b' := \sup (K \cap [\theta_-(H), b])$, and note
that $a\leq a' < b' \leq b$. Moreover, since $\psi_H$
is constant on each gap of $K$, we have $\psi_H(a')=\psi_H(a)$ and
$\psi_H(b')=\psi_H(b)$. Expand $a'$ and $b'$ in base 2 as
$$
a'=0. \theta_{t_1} \theta_{t_2} \cdots  \ \ \ \operatorname{and}
\ \ \ b'=0. \theta_{s_1} \theta_{s_2} \cdots,
$$
where $t_i,s_i \in \{ 0, 1 \}$ are uniquely determined by $a'$ and
$b'$. Set $t:=0.t_1 t_2 \cdots$ and $s:=0.s_1s_2 \cdots$ in base 2,
so that $t=\psi_H(a')$ and $s=\psi_H(b')$. Note that by the choice
of $a'$, if $t$ happens to be a dyadic rational, then the binary
expansion $0.t_1 t_2 \cdots$ is the one which terminates with a
string of $0$'s. Similarly, by the choice of $b'$, if $s$ is
dyadic, then the binary expansion $0.s_1 s_2 \cdots$ is the one
which terminates with a string of $1$'s. With this observation in
mind, let $j \geq 1$ be the smallest integer such that $t_j \neq
s_j$. Then,
\begin{equation}\label{eqn:up}
\psi_H(b)-\psi_H(a) = s-t \leq (\con) 2^{-j}.
\end{equation}
On the other hand, the open interval $]t,s[$ contains the dyadic point
$r=0.t_1\cdots t_{j-1} 1 \ov{0}$, so $]a',b'[$ contains the gap
$$
\left] 0.\theta_{t_1} \cdots \theta_{t_{j-1}} \theta_{0}
\ov{\theta_{1}}\, , \, 0.\theta_{t_1} \cdots \theta_{t_{j-1}}
\theta_{1} \ov{\theta_{0}} \right[
$$
whose endpoints are the two values of $A_H(r)$. Since the length of this gap
is bounded below by $(\con) 2^{-jp}$, it follows that
\begin{equation}\label{eqn:down}
b-a \geq b'-a' \geq (\con) 2^{-jp}.
\end{equation}
Combining the two inequalities (\ref{eqn:up}) and (\ref{eqn:down}), we obtain
$$\psi_H(b)-\psi_H(a) \leq (\con) (b-a)^{1/p},$$
which proves the lemma.
\end{pf}

\noindent
{\it Proof of \thmref{C}.} Let $F \subset {\cal R}$ denote the
(conjecturally empty) set of angles whose corresponding
parameter rays do not land. We prove that the set $A_H({\cal R} \sm F)
\subset K$ has Hausdorff dimension $1/p$. \thmref{C} would follow
since any other ray landing on $\eta_H$ which is not in $A_H({\cal R}
\sm F)$ must be rational of even denominator and there are countably
many such rays (compare the discussion before \lemref{Holder}).

Since $\dim_H({\cal R} \sm F)=1$ by \thmref{B},
Frostman's Lemma shows that for
any $0<\delta<1$ there exists a Borel probability measure $\mu$
supported in $\cal R \sm F$ such that $\mu(I) \leq (\con) \,
m(I)^\delta$ for all intervals $I$. Let $\nu:= (A_H)_{\ast} \mu$ be
the push-forward measure supported in $A_H({\cal R} \sm F)$. Take
any interval $J \subset [\theta_-(H), \theta_+(H)]$ and let
$I=\psi_H(J)$. Then $\nu(J)=\mu(I) \leq (\con) \, m(I)^{\delta}$.
By \lemref{Holder}, $m(I) \leq (\con) \, m(J)^{1/p}$. It follows
that $\nu(J) \leq (\con) \, m(J)^{\delta/p}$. By another
application of Frostman's Lemma, we conclude that
$$\frac{1}{p} = \Hdim(A_H(\TT)) \geq \Hdim(A_H({\cal R} \sm F))
\geq \frac{\delta}{p}\, .$$
Letting $\delta \to 1$, we obtain the result. \hfill $\Box$

\section{Critically non-recurrent real quadratics}
\label{sec:gen}

The results of the preceding section can be used to obtain, with
minimal effort, dimension estimates in the parameter space of real
quadratic polynomials. As an example, a difficult theorem of
Jakobson \cite{Jakobson} asserts that the nowhere dense set $\bd M
\cap \RR$ has positive linear measure, hence full Hausdorff
dimension 1. It is interesting to see that this last statement also
follows from \thmref{B}. To this end, let us recall the following
result from the theory of univalent maps (see \cite{Makarov} or
\cite{Pommerenke} for a proof):
\begin{quotation}
{\bf Makarov Dimension Theorem.} {\it For any univalent map $\phi :
\DD \to \CC$ and any Borel set $X \subset \TT$ with $\Hdim(X)=\delta$,
we have the estimates}
$$\Hdim (\phi(X)) > \left\{
\begin{array}{cl}
\ds{\frac{\delta}{2}} & \hspace{6mm} \operatorname{if\ \ }
\ds{0< \delta < \frac{11}{12}} \vspace{3mm}\\
\ds{\frac{\delta}{1+\sqrt{12(1-\delta)}}} & \hspace{6mm} \operatorname{if\ \ }
\ds{\frac{11}{12} < \delta < 1}
\end{array}
\right.$$
\end{quotation}
In our case, after removing a set $F$ of capacity zero
(conjecturally empty) from $\cal R$, we can assume that all the
parameter rays with angles in ${\cal R}\sm F$ land at a point of
$\bd M \cap \RR$. The set ${\cal R}\sm F$ still has Hausdorff
dimension 1, so by Makarov Dimension Theorem the image $\pi({\cal
R}\sm F) \subset \bd M \cap \RR$ has dimension at least 1. It
follows that $\Hdim (\bd M \cap \RR)=1$.

It was pointed out to me by S.~Smirnov that a more careful
application of the above argument, combined with the dimension
estimates in \S \ref{sec:AB}, shows the existence of a set of full
Hausdorff dimension in $\bd M \cap \RR$ which consists only of
critically non-recurrent quadratics. Recall that $c \in \CC$ is
called a {\it critically non-recurrent} parameter if the critical
point $0$ of the quadratic $Q_c$ does not belong to the closure of
its forward orbit $\{ Q_c^{\circ n}(0) \}_{n \geq 1}$. (Caution: In
the theory of interval maps, the term ``Misiurewicz'' is often used
for ``critically non-recurrent non-hyperbolic.'' Following the
standard terminology of complex dynamics, we have used the term
``Misiurewicz'' in the more restricted sense of ``critically finite
non-hyperbolic.'') It has been shown by D.~Sands that the set of
critically non-recurrent parameters in $\bd M \cap \RR$ has linear
measure zero \cite{Sands}. Here we prove a complement to his result
by showing that this set has
full Hausdorff dimension. \vspace{2mm}\\
{\it Proof of \thmref{CNR}.} Fix a small $\sigma >0$ such that
$(1-\sigma)/2 \notin {\cal R}$, and consider the compact set ${\cal
K}_{\sigma}$ defined in (\ref{eqn:EE}). By (\ref{eqn:rk}), ${\cal
K}_{\sigma} \cap [(1-\sigma)/2, 1/2] \subset {\cal R}$, so ${\cal
N}_{\sigma}:=\pi ({\cal K}_{\sigma} \cap [(1-\sigma)/2, 1/2])$ is a
well-defined subset of $\bd M \cap \RR$. We claim that every
parameter $c \in {\cal N}_{\sigma}$ is critically non-recurrent.
Assuming this is false, take a critically recurrent parameter $c$
in ${\cal N}_{\sigma}$ and a sequence $n_j \to \infty$ such that
$Q_c^{\circ n_j}(c) \to c$ as $j \to \infty$. By \thmref{2rays} the
critical value $c$ of $Q_c$ is the landing point of the two dynamic
rays at angles $\pm \tau(c)$, so $Q_c^{\circ n_j}(c)$ is the
landing point of the rays at angles $\pm \db^{\circ n_j}(\tau(c))$.
Since $Q_c^{\circ n_j}(c) \to c$, an easy exercise shows that
$\db^{\circ n_j}(\tau(c)) \to \pm \tau(c)$. In particular,
\begin{equation}\label{eqn:x}
\lim_{j \to \infty} \left|\db^{\circ n_j}(\tau(c))-\frac{1}{2} \right|
= \frac{1}{2} - \tau(c).
\end{equation}
On the other hand, the definition of ${\cal N}_{\sigma}$ and
the fact that $\tau =\pi^{-1}$ shows that
$\tau(c) \in {\cal K}_{\sigma} \cap [(1-\sigma)/2, 1/2]$. As in the proof of
\thmref{B}, it follows that
\begin{equation}\label{eqn:xx}
\left| \db^{\circ n}(\tau(c)) - \frac{1}{2} \right| \geq \frac{\sigma}{2}
\geq \frac{1}{2} - \tau(c)
\end{equation}
for all $n \geq 1$.
Comparing (\ref{eqn:x}) and (\ref{eqn:xx}), we obtain $\tau(c)=(1-\sigma)/2$.
This is a contradiction since $\tau(c) \in {\cal R}$ and $(1-\sigma)/2
\notin {\cal R}$.

Now choose a sequence $\sigma_n \to 0$ subject to the condition
$(1-\sigma_n)/2 \notin {\cal R}$; this is possible since $\cal R$
is nowhere dense. Use \lemref{E_sigma} and Makarov Dimension
Theorem to deduce $\lim_{n \to \infty} \Hdim {\cal N}_{\sigma_n}
=1$. This proves \thmref{CNR} since by the above argument the set
of critically non-recurrent parameters in $\bd M \cap \RR$ contains
${\cal N}_{\sigma_n}$ for all $n$. \hfill $\Box$

\section{Biaccessibility in real quadratics}
\label{sec:bi}

Let $P:\CC \to \CC$ be a polynomial with connected Julia set $J(P)$.
A point $z
\in J(P)$ is called {\it biaccessible} if there are two or more
dynamic rays landing at $z$. We denote the set of all such points
by $\tilde{J}(P)$. It is known that the harmonic measure of
$\tilde{J}(P)$ is zero unless $P$ is affinely conjugate to a
Chebyshev polynomial for which $J(P)$ is an interval and hence
$\tilde{J}(P)$ has full harmonic measure $1$. This was proved by
the author for all locally-connected and some non locally-connected
quadratic Julia sets \cite{Zakeri}. Proofs for the general case
were later found independently by S.~Smirnov \cite{Stas} and
A.~Zdunik \cite{Zdunik}.

Naturally, one would like to know about the Hausdorff dimension of
the set of angles landing at biaccessible points of a non-Chebyshev
polynomial. In this section we prove \thmref{D} which shows how this
dimension can be effectively estimated, at least for real quadratics.
Let $c \in [-2,1/4]$, $J_c=J(Q_c)$,
$\tilde{J}_c=\tilde{J}(Q_c)$, and let $\beta > 0$ be the fixed
point of $Q_c$ at which the dynamic ray $\rd(0)$ lands. Define
$$\begin{array}{rl}
S_c  := & \{ t \in \TT : \rd(t) \ \text{lands on} \ [-\beta, \beta]
\} , \vspace{1mm} \\
B_c := & \{ t \in \TT : \rd(t) \ \text{lands at a point of} \ \tilde{J}_c \}.
\end{array}$$
It is not hard to prove that
$$\tilde{J}_c=\bigcup_{n \geq 0} Q_c^{-n}( J_c \, \cap \,
]-\beta, \beta[).$$
(see \cite{Zakeri}), so that
\begin{equation}
\label{eqn:equal0}
\Hdim (\tilde{J}_c) = \Hdim (J_c \, \cap \,  ]-\beta, \beta[).
\end{equation}
By ignoring the countable set of angles of
dynamic rays landing at $\beta$ and its iterated preimages, it also
follows that
\begin{equation}
\label{eqn:equal}
\Hdim (B_c) = \Hdim (S_c).
\end{equation}
In order to prove \thmref{D}, we link the set $S_c$ to the family of
compact sets $\{ {\cal K}_{\sigma} \}_{\sigma >0}$ introduced in \S \ref{sec:AB}.
Once this connection is found, we simply use the dimension
estimates in \S \ref{sec:AB} to find bounds on the dimension of
$S_c$.

\begin{lemma}[Douady]
\label{doo}
If $c \in \bd M \cap \RR$, then
\begin{equation}\label{eqn:S_c}
S_c=\{ t \in \TT : \db^{\circ n}(t) \notin \, ] \tau(c) , 1-\tau(c) [ \
\operatorname{for \ all \ } n \geq 1 \}.
\end{equation}
\end{lemma}

\begin{pf}
We follow the argument in \cite{Douady3}. Let $E$ denote the closed set
defined by the right
side of (\ref{eqn:S_c}). If $t \in S_c$, then $\rd(t)$ lands at
some $z \in [-\beta, \beta]$, and so for $n \geq 1$ the dynamic ray
$\rd(\db^{\circ n}(t))$ lands at $Q_c^{\circ n}(z)$. Evidently,
$Q_c^{\circ n}(z) \geq c$, and since $J_c \, \cap \, ]c,r_c[ =
\emptyset$, we actually have the stronger inequality $Q_c^{\circ
n}(z) \geq r_c$. Since $r_c$ is the landing point of the rays $\rd(\pm \tau(c))$,
it follows that $\db^{\circ n}(t) \notin \, ]
\tau(c) , 1-\tau(c) [$ for all $n \geq 1$, which implies $t \in E$.
This proves $S_c \subset E$. On the other hand, let $]s,t[$ be a
connected component of $\TT \sm E$. Then there exists an $n \geq 1$
so that the iterate $\db^{\circ n}$ maps $]s,t[$ homeomorphically
to $]\tau(c), 1-\tau(c)[$, but $\db^{\circ j}(]s,t[) \, \not
\subset \, ]\tau(c), 1-\tau(c)[$ if $1 \leq j \leq n-1$. It easily
follows by a recursive argument that for every $0 \leq j \leq n$
the dynamic rays $\rd(\db^{\circ j}(s))$ and $\rd(\db^{\circ
j}(t))$ land on $[-\beta, \beta]$. In particular, $s,t \in S_c$.
Since $E$ has no interior, the set of all such $s,t$ is dense in
$E$, proving $E \subset S_c$.
\end{pf}

\begin{corollary}\label{S=K}
For every $c \in \bd M \cap \RR$, we have $S_c={\cal K}_{\sigma}$,
where $\sigma = 1-2\tau(c)$ and where ${\cal K}_{\sigma}$ is the
compact set defined in (\ref{eqn:EE}).
\end{corollary}

\begin{pf}
This simply follows from the previous lemma and the definition of
${\cal K}_{\sigma}$.
\end{pf}

\begin{corollary}\label{es1}
For $c \in \bd M \cap [-2,-1.75]$, we have
$$\Hdim (B_c) \geq 1 - \frac{1}{\log 2} \, \log \left( \frac{16\tau(c)-5}
{32\tau(c)-13}\right) $$
\end{corollary}

Observe that $c=-1.75$ is the landing point of the
parameter ray $\rp(3/7)$, i.e., the root point of the unique
real hyperbolic component of period $3$.

\begin{pf}
If $c=-2$, then $\tau(c)=1/2$ and the inequality is evident. So
assume $3/7 \leq \tau(c) < 1/2$, and choose
a dyadic $\sigma = 2^{-p} \leq 1/4$ such that
$\sigma /2 < 1-2 \tau(c) \leq \sigma$. Then, by
\corref{S=K}, $S_c \supset {\cal K}_{\sigma}$. By \lemref{aux},
\begin{align*}
\Hdim(S_c) \geq \Hdim({\cal K}_{\sigma}) & \geq 1-\frac{1}{\log 2}
\, \log \left( \frac{3-4 \sigma}{3-8 \sigma} \right) \vs \\
& \geq 1 - \frac{1}{\log 2} \, \log \left( \frac{16\tau(c)-5}
{32\tau(c)-13}\right).
\end{align*}
Since $\Hdim (S_c) =\Hdim (B_c)$ by (\ref{eqn:equal}), we obtain
the result.
\end{pf}

We can improve the above estimate and at the same time extend it to
all real parameters in $[-2,-1.75]$. To this end, we need a few
preliminary remarks.

Recall that ${\cal H}$ is the collection of all real hyperbolic
components of $M$. Let $H \in {\cal H}$ and let $H_1$ be the
$1/2$-satellite of $H$ (see \S \ref{sec:bg} for the definition).
Then $H_1 \in {\cal H}$ and in the passage from $H$ to $H_1$ the
attracting cycle of $Q_c$ undergoes a period-doubling bifurcation.
We express the relationship between $H_1$ and $H$ by writing $H_1
\til H$.

\begin{lemma}\label{propp}
Let $H, H_1 \in {\cal H}$ have root points $c$ and $c_1$, and
assume $H_1 \til H$. Consider any real parameter $c_0 \in H$. Then
\begin{enumerate}
\item[(i)]
$S_{c_0}=S_c$ and $B_{c_0}=B_c$.
\item[(ii)]
$S_{c_1} \sm S_c$ is countably infinite.
\end{enumerate}
\end{lemma}

\begin{pf}
Statement (i) follows from the well-known fact that the quadratics
$Q_c$ and $Q_{c_0}$ are combinatorially equivalent, so that
$\rd(t)$ and $\rd(s)$ land at a common point in the dynamic plane
of $Q_c$ if and only if $R_{c_0}(t)$ and $R_{c_0}(s)$ land at a
common point in the dynamic plane of $Q_{c_0}$ (see for example
\cite{Milnor}). To see (ii), first note $\tau(c_1) > \tau(c)$ so
that $S_{c_1} \supset S_{c}$ by \lemref{doo}. Also note that
$c_1=\iota_H(-0.75)$, where $-0.75$ is the root point of the
$1/2$-satellite of the main cardioid. Hence, by the discussion of
tuning in \S \ref{sec:bg}, the filled Julia set $K_{c_1}$ is
obtained from $K_{c_0}$ by replacing each bounded Fatou component
of $Q_{c_0}$ by a copy of the filled Julia set $K_{-0.75}$. A brief
inspection shows that
$$
S_{-0.75}= \{ 0 , 1/2 \} \cup \left\{ \pm \frac{1}{6\cdot 2^n}
\right\} _{n\geq 0} \cup \left\{ \frac{1}{2} \pm \frac{1}{6\cdot
2^n} \right\} _{n\geq 0},
$$
so the set $S_{-0.75}$ is countably infinite. Thus, the process of
constructing $K_{c_1}$ from $K_{c_0}$ adds only countably many new
angles to $S_{c_0}=S_c$, and we obtain the result. The details are
straightforward and will be left to the reader.
\end{pf}

From \lemref{doo} and \lemref{propp}(i) one immediately obtains the
monotonicity of the family $S_c$. In particular,

\begin{corollary}
The dimension function $c \mapsto \Hdim(B_c)=\Hdim(S_c)$ is
monotonically decreasing on $[-2,1/4]$.
\end{corollary}

Given any $H \in {\cal H}$, one can form the infinite cascade
\begin{equation}
\label{eqn:cas}
\cdots \til H_n \til \cdots \til H_2 \til H_1 \til H
\end{equation}
of period-doubling bifurcation components originating from $H$.
Then the root $c_n$ of $H_n$ converges exponentially fast to a
limit $c^*(H) \in \bd M \cap \RR$ as $n \to \infty$ (compare
\cite{Lyubich2}).

\begin{definition}
We call the parameter $c^*(H)$ {\it the main Feigenbaum point
associated to} $H$.
\end{definition}

As an example, the much studied main Feigenbaum point associated to
the main cardiod $H_0$, which we denote by
$c_{\operatorname{Feig}}$, is approximately located at $-1.401155$
(compare \figref{Feig}).

It is not hard to see using the properties
of the tuning map that for every $H \in {\cal H}$,
$$c^*(H)= \iota_{H}(c_{\operatorname{Feig}}).$$
Note that the assignment $H \mapsto c^*(H)$ is far from being
one-to-one. In fact, $c^*(H)=c^*(H_1)$ whenever $H_1 \til H$, and so
all the components in (\ref{eqn:cas}) have the same main Feigenbaum
point associated to them.

Using the fact that $S_{1/4}=\{ 0, 1/2 \}$ and applying
\lemref{propp} repeatedly, we obtain

\begin{corollary}
For $c_{\operatorname{Feig}} < c \leq 1/4$, $\Hdim
(B_c)=\Hdim(S_c)=0$.
\end{corollary}

The following definition will be used in the proof of \thmref{D}.

\begin{definition}
For $c \in [-2,1/4]$, we define
$$\rho(c):= \left\{ \begin{array}{ll}
\tau(c^*(H)) & \text{if} \ c \in \ov{H}
\ \text{for some}\ H \in {\cal H} \vs \\
\tau(c) & \text{otherwise}
\end{array}
\right.$$
In other words, if $c$ is not hyperbolic or parabolic,
then $\rho(c)=\tau(c)$, but if $c$ is hyperbolic or parabolic,
then $\rho(c)$ is the external angle of the associated main
Feigenbaum point. Clearly $\rho(c) \geq \tau(c)$ whenever
$c \in \bd M \cap \RR$.
\end{definition}

We are now ready to prove \thmref{D} cited in the introduction. \vs \\
{\it Proof of \thmref{D}.} For $-2 < c \leq -1.75$, define
\begin{equation}\label{eqn:ell}
\ell(c):=1 - \frac{1}{\log 2} \log \left( \frac{16\rho(c)-5}
{32\rho(c)-13} \right).
\end{equation}
If $c$ does not belong to the closure of any real hyperbolic
component, then $\ell(c)$ coincides with the lower bound in
\corref{es1}, so $\Hdim(B_c) \geq \ell(c)$. If, on the other hand, $c
\in \ov{H}$ for some $H \in {\cal H}$, then we may consider the
cascade of period doubling bifurcations (\ref{eqn:cas}), with the
root points $c_n$. An inductive application of \lemref{propp} then
shows that $S_{c_n} \sm S_c$ is countably infinite for every $n$.
Hence, by \corref{es1},
$$\Hdim(B_c)=\Hdim(S_c)=\Hdim(S_{c_n})=\Hdim(B_{c_n}) \geq
1 - \frac{1}{\log 2} \, \log \left( \frac{16\tau(c_n)-5} {32\tau(c_n)-13}
\right).$$
But $c_n
\to c^*(H)$ implies $\tau (c_n) \to \tau (c^*(H))=\rho(c)$ by
\lemref{R=c}, so that the right side of the above inequality tends
to $\ell(c)$ as $n \to \infty$, implying $\Hdim(B_c) \geq \ell(c)$. By
\lemref{<1} and \corref{S=K}, $\Hdim(S_c)<1$. Combining these two
inequalities, we obtain
$$\ell(c) \leq \Hdim(B_c) <1.$$
The fact that $-2 < c \leq -1.75$ implies $4/9 < \rho(c) < 1/2$,
and a brief computation shows that $\ell(c)>0$. That
$\lim_{c\,{\scriptstyle\searrow}\, -2} \ell(c) =1$ follows from
$\lim_{c\,{\scriptstyle\searrow}\, -2}\rho(c)= 1/2$. \hfill $\Box$

\begin{remark}
I do not know if $\Hdim (B_{c_{\operatorname{Feig}}}) >0$. More
generally, I do not know how to obtain explicit lower bounds for
$\Hdim (B_c)$ for $-1.75 < c \leq c_{\operatorname{Feig}}$. On the
other hand, the preceding arguments show $\Hdim(B_c)<1$ for $-2 < c
\leq 1/4$ and I do not know if this phenomenon is general in the
quadratic family. In other words, does there exist a complex
parameter $c \neq -2$ for which $\Hdim (B_c)=1$?
\end{remark}

Finally, it is quite easy to prove an analogue of \thmref{D} for
the set of biaccessible points in the Julia set itself.

\begin{corollary}
\label{HDJ}
Let $-2 < c \leq 1/4$ and let $\tilde{J}_c$ denote the set of
biaccessible points in the Julia set $J_c$. If $J_c$ is full, then
$\Hdim (\tilde{J}_c)=1$. On the other hand, if $J_c$ is not full
and $c \leq -1.75$, then
\begin{equation}
\label{eqn:3}
0<\ell'(c) \leq \Hdim(\tilde{J}_c) \leq 1.
\end{equation}
where $\ell'(c)$ is an explicit constant which tends to $1$ as
$c$ tends to $-2$. In particular,
\begin{equation}
\label{eqn:4}
\lim_{c\,{\scriptstyle\searrow}\, -2} \Hdim (\tilde{J}_c) = \Hdim
(\tilde{J}_{-2}) =1.
\end{equation}
\end{corollary}

\begin{pf}
If $J_c$ is full, then $J_c \supset [ -\beta, \beta ]$ and hence
$\Hdim (\tilde{J}_c)=1$ by (\ref{eqn:equal0}). If $J_c$ is not full
and $-2 < c \leq -1.75$, then we can take $\ell'(c):=\lambda
(\ell(c))$, where $\ell(c)$ is given by (\ref{eqn:ell}) and
$\lambda=\lambda(\delta)$ is the lower bound function given by
Makarov Dimension Theorem. The upper bound in (\ref{eqn:3}) follows
from (\ref{eqn:equal0}). Finally,
$\lim_{c\,{\scriptstyle\searrow}\, -2} \Hdim (B_c) =1$ together
with Makarov Dimension Theorem shows that
$\liminf_{c\,{\scriptstyle\searrow}\, -2} \Hdim (\tilde{J}_c) \geq
1$, which combined with $\Hdim (\tilde{J}_c) \leq 1$ proves
(\ref{eqn:4}).
\end{pf}

\end{document}